\documentclass[fleqn,a4paper,10pt]{amsart}

\usepackage{graphicx}
\usepackage{tikz}
\usetikzlibrary{shapes,decorations.pathmorphing}
\usepackage{amsmath}
\usepackage{amssymb,amsthm}
\usepackage{comment}
\usepackage{hyperref}

\theoremstyle{definition}

\newtheorem{proper}{Property}

\title[Two-point function of quadrangulations:
a new derivation]{The distance-dependent two-point function of quadrangulations:
a new derivation by direct recursion}
\author{Emmanuel Guitter}
\address{Institut de Physique Th\'eorique, CEA, IPhT, 91191 Gif-sur-Yvette, France, CNRS, UMR 3681}
\email{emmanuel.guitter@cea.fr}

\begin{document}
\maketitle

\begin{abstract}
We give a new derivation of the distance-dependent two-point function of planar quadrangulations by 
solving a new direct recursion relation for the associated slice generating functions. Our approach for both the derivation 
and the solution of this new recursion is in all points similar to that used recently by the author in the context 
of planar triangulations.  
\end{abstract}

\section{Introduction}
\label{sec:introduction}
The \emph{distance-dependent two-point function} of a family of maps is, so to say, the generating function of these maps
with two marked ``points" (e.g.\ vertices or edges) at a prescribed graph distance from each other. It informs us about the 
distance profile between pairs of points picked at random on a random map in the ensemble at hand.
In the case of planar maps, 
explicit expressions for the distance-dependent two-point function of a number of map families were obtained by several techniques 
\cite{GEOD,PDFRaman,BG12,AmBudd,BFG,FG14},
all based on the relationship which exists between the two-point function and generating functions 
for either some particular decorated trees, or equivalently for some particular pieces of maps called \emph{slices}.  
This relationship is itself a consequence of the existence of some now well-understood bijections between maps and trees or slices
\cite{SchPhD,BDG04}.

In a recent paper \cite{G15}, we revisited the distance-dependent two-point function of planar triangulations (maps whose all 
faces have degree $3$) and showed how to obtain its expression from the solution of some direct recursion relation 
on the associated slice generating functions. The solution of the recursion made a crucial use of some old results by Tutte in his seminal
paper \cite{TutteCPT} on triangulations. In this paper, we extend the analysis of \cite{G15} to the case of planar quadrangulations
(maps which all faces of degree $4$) by showing that a similar recursion may be written and solved by the same treatment as for triangulations. 

The paper is organized as follows: we start in Section~\ref{sec:slices} by giving the basic definitions (Sect.~\ref{sec:basics})
and by recalling
the relation which exists between the distance-dependent two-point function of planar quadrangulations and the 
generating functions of particular slices (Sect.~\ref{sec:GkRk}). We then derive in Section~\ref{sec:recur} a direct recursion relation for 
the slice generating functions, based on the definition of a particular \emph{dividing line} drawn on the slices (Sect.~\ref{sec:dividing})
and on a decomposition of the slices along this line (Sect.~\ref{sec:decomp}). 
Section~\ref{sec:simplemaps} shows how to slightly simplify the recursion by reducing the problem to slice generating 
functions for \emph{simple quadrangulations}, i.e.\ quadrangulations without multiple edges (Sect.\ref{sec:substitution}). This allows to 
make the recursion relation fully tractable by giving an explicit expression for its kernel (Sects.~\ref{sec:eqforPhi} and \ref{sec:kernel}). 
Section~\ref{sec:final} 
is devoted to solving the recursion relation, first in the case of simple quadrangulations (Sect.~\ref{sec:simplesol}), then for general ones
(Sect.~\ref{sec:general}),
leading eventually to some explicit expression for the distance-dependent two-point function. We conclude in Section~\ref{sec:conclusion}
with some final remarks.

\section{The two-point function and slice generating functions}
\label{sec:slices}
\subsection{Basic definitions}
\label{sec:basics}
As announced, the aim of this paper is to compute the distance-dependent two-point function of planar quadrangulations.
Recall that a planar map is a connected graph embedded on the sphere. The map is \emph{pointed} if it has a marked 
vertex (the pointed vertex) and  \emph{rooted} if it has a marked oriented edge (the root-edge). In this latter case, the origin 
of the root-edge is called the root-vertex. A planar quadrangulation is a planar map whose all faces have degree $4$. 
For $k\geq 1$, we define the distance-dependent two-point function $G_k\equiv G_k(g)$ of planar quadrangulations  
as the generating function of pointed rooted quadrangulations whose pointed vertex and root-vertex are at graph distance 
$k$ from each other. The quadrangulations are enumerated with a weight $g$ per face.
Note that, since planar quadrangulations are bipartite maps, the graph distances from a given vertex
to two neighboring vertices have different parities, hence their difference is $\pm1$. In particular, in quadrangulations 
enumerated by $G_k$, the endpoint of the root-edge is necessarily at distance $k-1$ or $k+1$ from the pointed vertex.

A \emph{quadrangulation with a boundary} is a rooted planar map whose all faces have degree
$4$, except the root-face, which is the face lying on the right of the root-edge, which has arbitrary degree.
Note that this degree is necessarily even as the map is clearly bipartite. 
The faces different from the root-face are called \emph{inner faces} and form the \emph{bulk} of the map while the edges incident to the root-face 
(visited, say clockwise around the bulk) form the \emph{boundary} of the map, whose \emph{length} is the degree of the root-face.
\begin{figure}
\begin{center}
\includegraphics[width=6cm]{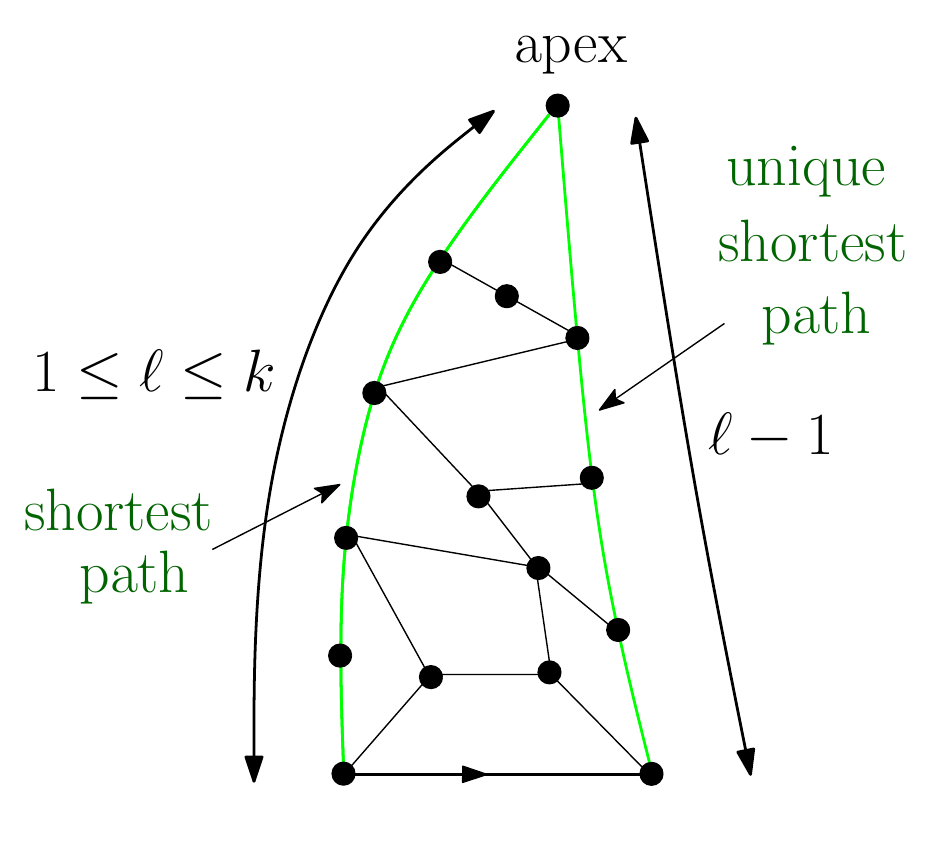}
\end{center}
\caption{A example of slice with a boundary of length $2\ell=10$, hence with left-boundary length $\ell=5$, thus
contributing to $R_k$ for all $k\geq 5$.}
\label{fig:slice}
\end{figure}

As in \cite{G15}, we may compute $G_k$ by relating it to the generating function of \emph{slices}, which 
are particular instances of quadrangulations with a boundary, characterized by the following properties:
let  $2\ell$ ($\ell\geq 1$) be the length of the boundary, we call \emph{apex} the vertex reached from the root-vertex 
by making $\ell$ elementary steps \emph{along the boundary} clockwise around the bulk. The map at hand is a slice if
(see figure \ref{fig:slice}):
 \begin{itemize}
 \item the graph distance from the root-vertex to the apex is $\ell$.
Otherwise stated, the left boundary of the slice, which is the portion (of length $\ell$) of boundary between the root-vertex and the apex
 clockwise around the bulk  is a shortest path between its endpoints within the map;
 \item the distance from the endpoint of the root-edge to the apex is $\ell-1$.
Otherwise stated, the right boundary of the slice, which is the portion (of length $\ell-1$) of boundary between the endpoint of 
the root-edge and the apex 
 counterclockwise around the bulk is a shortest path between its endpoints within the map; 
 \item the right boundary is the \emph{unique} shortest path between its endpoints within the map;
 \item the left and right boundaries do not meet before reaching the apex.
 \end{itemize}
 \begin{figure}
\begin{center}
\includegraphics[width=4cm]{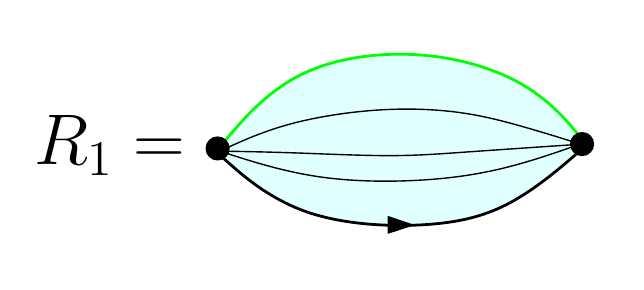}
\end{center}
\caption{A schematic picture of a map enumerated by $R_1$, referred to in this paper as a \emph{bundle} between its two
boundary vertices. We indicated only those edges which connect the extremities of the root-edge to emphasize the fact that 
their number may be arbitrarily large.}
\label{fig:R1}
\end{figure}
We call $R_k\equiv R_k(g)$ ($k\geq 1$) the generating function of slices with $1\leq \ell \leq k$, enumerated with a weight $g$
per inner face. Note that the root-edge-map, which is the map reduced to the single root-edge and a root-face of degree $2$ is a slice with $\ell=1$ and
contributes a term $1$ to all $R_k$ for $k\geq 1$.
The generating function $R_1$ deserves some special attention: by definition, $R_1$ enumerates
slices with $\ell=1$, hence with a boundary of length $2$. The right boundary has length $0$ and the apex is
the endpoint of the root edge while the left boundary, of length $1$, connects  both extremities of the root-edge (which are 
necessarily distinct). This connection is also performed by the root-edge itself, and the map forms in general what we shall call a a
 \emph{bundle} between the extremities of the root-edge (see figure~\ref{fig:R1}). 
The function $R_1$ is thus the generating function of bundles between adjacent vertices.
\vskip .2cm
\begin{figure}
\begin{center}
\includegraphics[width=8cm]{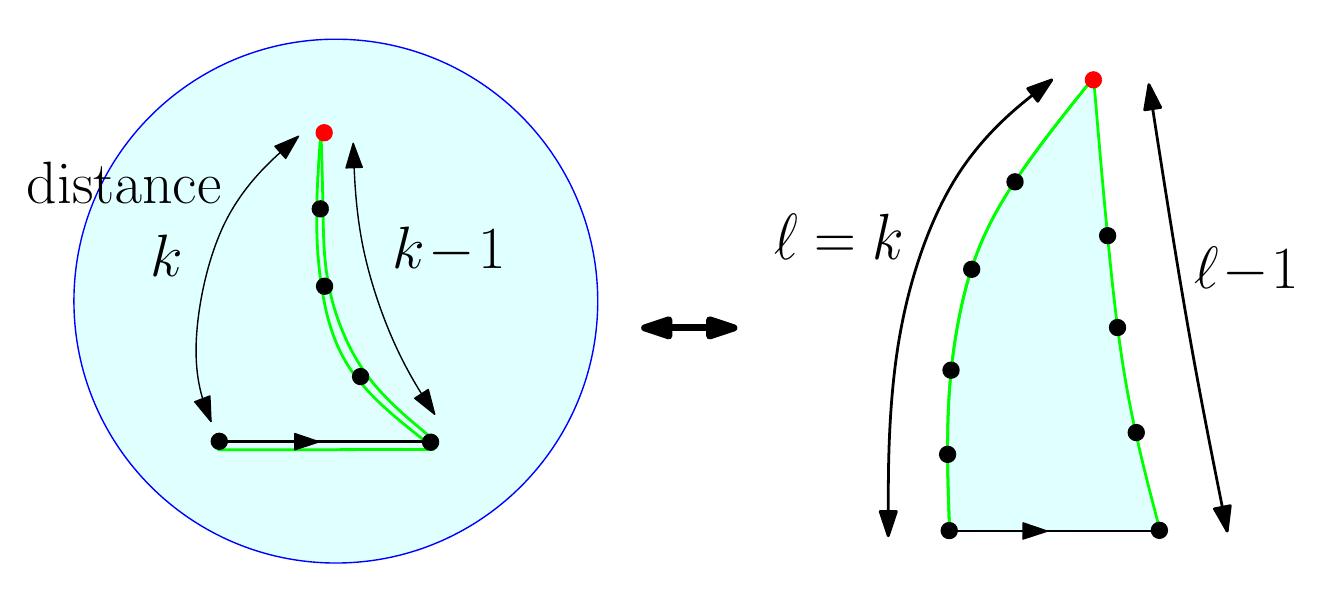}
\end{center}
\caption{The one-to-one correspondence between pointed rooted quadrangulations whose root-edge 
has its extremities at respective distance $k$ and $k-1$ from the pointed vertex (in red) and a slice with left-boundary length $\ell=k$.
On the left, we have drawn (in green) the leftmost shortest paths (starting with the root-edge itself) from the the root-vertex to the pointed vertex. 
Cutting along this path builds the slice on the right.}
\label{fig:twopoint}
\end{figure}

\subsection{Relation between $\mathbf{G_k}$ and $\mathbf{R_k}$}
\label{sec:GkRk}
We may now easily relate $G_k$ to $R_k$ via the following argument: consider a pointed rooted quadrangulation enumerated by $G_k$.
As already mentioned, the endpoint of the root-edge is necessarily at distance $k-1$ or $k+1$ from the pointed vertex, which divides the
maps at hand into two categories. Assume that the map belongs to the first category, for which the distance is $k-1$. 
Then we may draw the \emph{leftmost} shortest path from the root-vertex to the pointed vertex, choosing
as first step the root-edge itself (see figure~\ref{fig:twopoint}). Cutting along this shortest path creates a map with a boundary of length $2k$ which
is easily seen to be a slice of left-boundary length $k$, which is moreover not reduced to the root-edge-map when $k=1$. 
Such slices are enumerated by $R_k-R_{k-1}$ (since we must suppress from $R_k$ the slices with $1\leq \ell\leq k-1$) for $k\geq 2$ and by $R_1-1$ for $k=1$. This yield a contribution $R_k-R_{k-1}-\delta_{k,1}$ to $G_k$ from the first category, where we take the usual convention that $R_0=0$.
The second category corresponds to maps whose root-edge has its endpoint at distance $k+1$ from the pointed vertex.
By reversing the orientation of the root-edge, they are in bijection with maps of the first category, up to a change $k\to k+1$, hence
are enumerated by $R_{k+1}-R_{k}$ (since $\delta_{k+1,1}=0$ for $k\geq 1$). To summarize, we have the relation
\begin{equation}
G_k=(R_k-R_{k-1}-\delta_{k,1})+(R_{k+1}-R_{k})=R_{k+1}-R_{k-1}-\delta_{k,1}\ , \qquad k\geq 1
\label{eq:GkRk}
\end{equation}
which the convention $R_0=0$.

As for the slice generating functions $R_k$ ($k\geq 1$), they satisfy the now well-known equation \cite{GEOD} (see below for its derivation)
\begin{equation}
R_k=1+g\, R_k(R_{k-1}+R_k+R_{k+1})\ , \qquad k\geq 1
\label{eq:Rkeq}
\end{equation}   
with $R_0=0$. In particular, it is interesting to introduce the quantity $R_\infty=\lim_{k\to\infty} R_k$ which is the 
generating function of slices with arbitrary left-boundary length $\ell\geq 1$. From \eqref{eq:Rkeq}, this quantity is
directly obtained as the solution of
\begin{equation}
R_\infty=1+3g\, R_\infty^2
\label{eq:Rinf}
\end{equation}  
which satisfies $R_\infty=1+O(g)$. Equation \eqref{eq:Rkeq} may be viewed as a recursion on $k$ (giving $R_{k+1}$ from the knowledge 
of $R_k$ and $R_{k-1}$) but this recursion requires at initial data the knowledge of $R_1$. In the present case, it can
be shown \cite{CENSUS} that $R_1=R_\infty-g\, R_\infty^3$ and \eqref{eq:Rkeq} allows one in principle to determine $R_k$ for all $k$. Getting 
an \emph{explicit expression} for $R_k$ by this approach is a different story and so far, no real constructive way to 
solve \eqref{eq:Rkeq} was proposed. Instead, the method used so far was to first \emph{guess} the expression for $R_k$ and then 
verify that it solves \eqref{eq:Rkeq}. This led to the explicit formula for $R_k$ given in \cite{GEOD}, and eventually to $G_k$.
Another approach to determine $R_k$ was elaborated in \cite{BG12} where it was shown that the $R_k$'s appear as coefficients 
in a suitable continued fraction expansion for a standard generating function of quadrangulations with a boundary.
In this paper, we present a new direct recursion relation for $R_k$ which we shall then solve explicitly in a constructive way.
\begin{figure}
\begin{center}
\includegraphics[width=11cm]{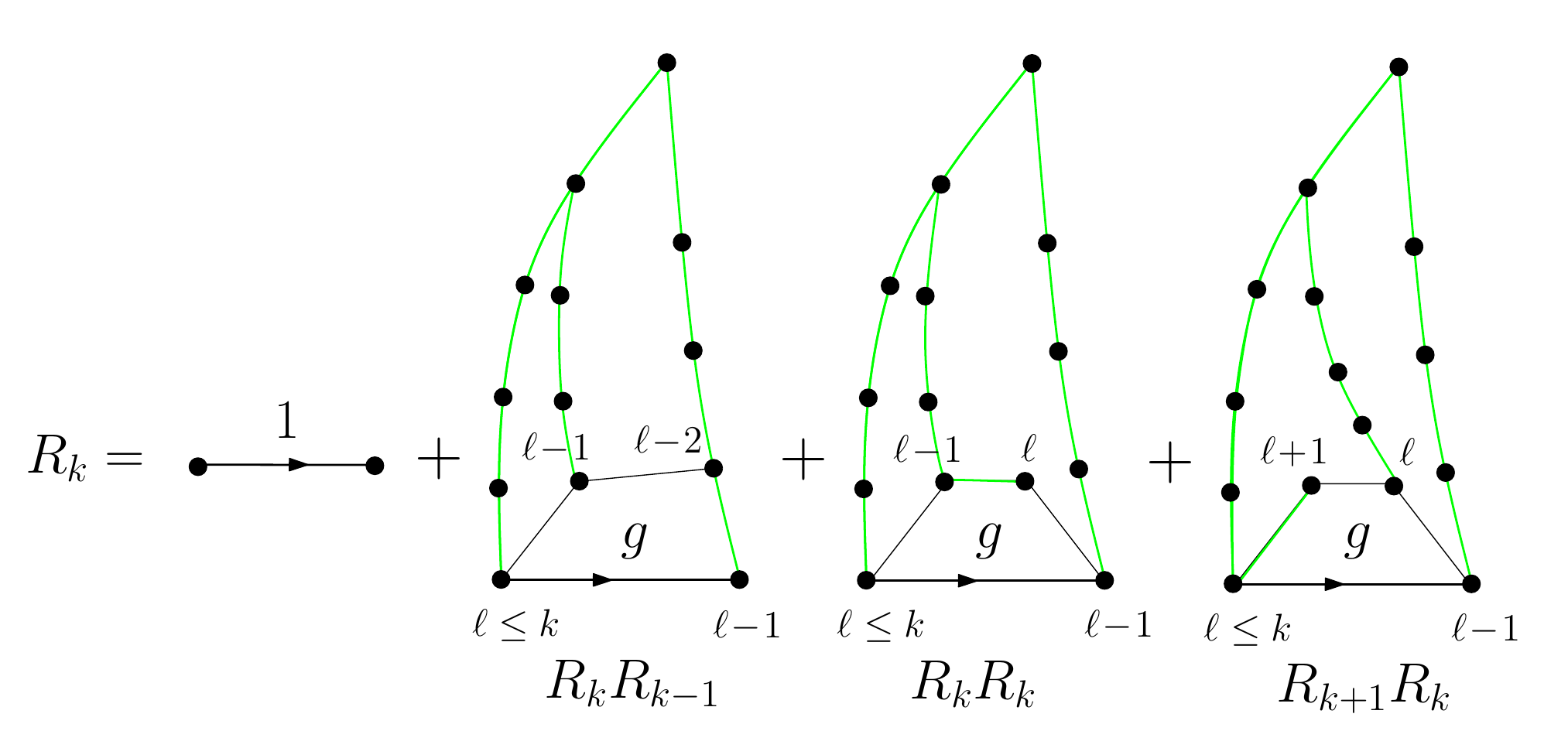}
\end{center}
\caption{The decomposition of slices leading to eq.~\eqref{eq:Rkeq}. If not reduced to the root-edge-map (weight $1$), the slice with $1\leq\ell\leq k$ 
is decomposed by removing the face immediately on the left of the root-edge (weight $g$), whose intermediate vertices are at distance $(\ell-1,\ell-2)$
(in which case the second vertex lies on the right boundary),
$(\ell-1,\ell)$ or $(\ell+1,\ell)$. 
Cutting along the leftmost shortest paths from these vertices to the apex produces two sub-slices
enumerated by $R_kR_{k-1}$, $R_k^2$ and $R_kR_{k+1}$ respectively.}
\label{fig:recur}
\end{figure}

To end this section, let us briefly recall for completeness the derivation of \eqref{eq:Rkeq}.
Consider a map enumerated by $R_k$ not reduced to the root-edge-map and consider the face directly on the left of the
root-edge. If the left-boundary length of the slice is $\ell$ ($1\leq \ell \leq k$), the sequence of distances to the apex of the 
four successive vertices\footnote{In all generality, it may happen that the four vertices are not distinct, in which case we must more precisely
consider the four successive \emph{corners} around the face, the distance to the apex of a corner being the distance to the apex of the incident vertex. 
Our statements may be straightforwardly adapted to these cases.} of the face, clockwise around the face starting from the root-vertex is either $\ell\to\ell-1\to\ell-2\to\ell-1$ (if $\ell\geq 2$), $\ell\to\ell-1\to\ell\to\ell-1$ or $\ell\to\ell+1\to\ell\to\ell-1$ (see figure~\ref{fig:recur}). 
Note that each of these paths of length $3$ has exactly two ``descending
steps" (i.e.\ steps for which the distance decreases by $1$). We may now draw, starting from the two intermediate vertices, the leftmost shortest paths from 
these vertices to the apex. This divides the slice into two slices (see figure~\ref{fig:recur}) whose two root-edges correspond precisely to the two descending steps.
For the sequence $\ell\to\ell-1\to\ell-2\to\ell-1$, the respective left-boundary lengths of the two slices are $\ell'$ and $\ell''-1$ with $\max(\ell',\ell'')=\ell$.
Demanding $\ell\leq k$ is equivalent to demanding $\ell'\leq k $ and $\ell''-1\leq k-1$ so that the slice pairs are enumerated
by $R_k\, R_{k-1}$, which explains the first of the three quadratic terms in \eqref{eq:Rkeq}. The two other terms come from the
two other possible distance sequences, while the first term $1$ corresponds to the edge-root-map. This explains \eqref{eq:Rkeq}.
  

\section{A direct recursion relation for slice generating functions}
\label{sec:recur}
\subsection{Definition of the dividing line}
\label{sec:dividing}
\begin{figure}
\begin{center}
\includegraphics[width=7cm]{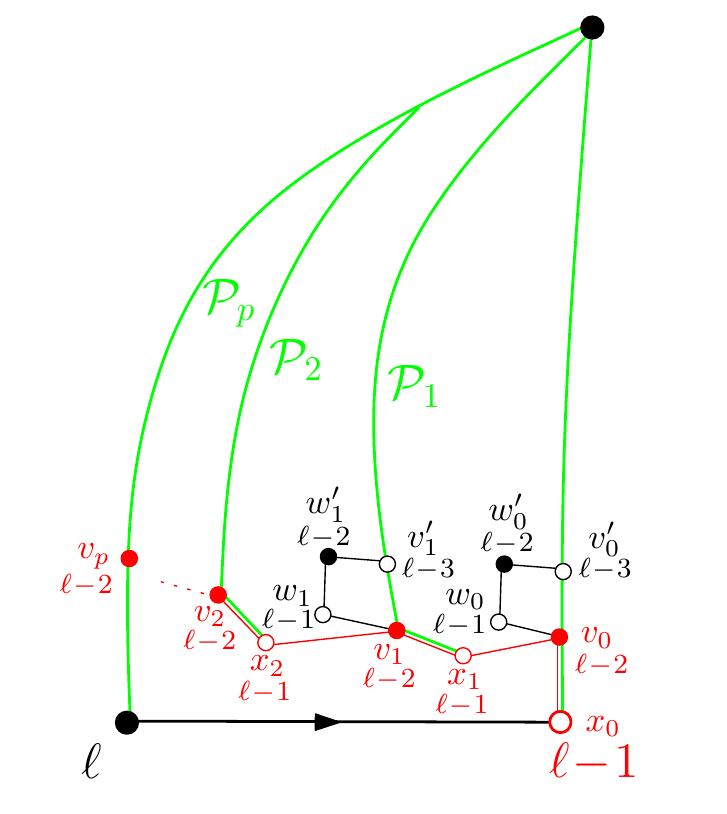}
\end{center}
\caption{The construction of the dividing line (see text for details).}
\label{fig:dividing}
\end{figure}
We shall now derive a new direct recursion for $R_k$. More precisely, our recursion is best expressed in terms of the generating
function 
\begin{equation}
T_k\equiv R_k-R_1\ , \qquad k\geq 1
\end{equation}
which enumerates slices with left-boundary length $\ell$ satisfying $2\leq \ell \leq k$ (in particular $T_1=0$). 
As in \cite{G15}, our recursion is based on a decomposition of 
the map along a particular \emph{dividing line} which we define now. Consider a slice with left-boundary length $\ell\geq 3$.
The dividing line is a sequence of edges which are alternatively of type $\ell-2\to \ell-1$ and $\ell-1\to \ell-2$. It forms 
a simple open curve which connects the right and left boundaries of the slice, hence separates the apex from the root-vertex. 
To construct the line, we proceed as follows: consider the vertices $v_0$ and $v_0'$ of the right boundary at respective distance 
$\ell-2$ and $\ell-3$ from the apex and consider the face directly on the left of the right-boundary edge linking $v_0$ to $v_0'$
(see figure~\ref{fig:dividing}).
Call $w_0$ and $w_0'$ the two other vertices incident to this face (so that the sequence clockwise around the face is $v_0 \to w_0 \to w_0'\to v_0'$).
The four vertices are necessarily distinct as otherwise, we could find a shortest path from $v_0$ to the apex lying strictly to the left of 
the right boundary, in contradiction with the fact that the right-boundary is the unique shortest path from the endpoint $x_0$
of the root-edge to the apex. Moreover, the distance from $w_0$ to the apex (which is {\it a priori} $\ell-3$ or $\ell-1$) cannot be equal
to $\ell-3$ as this would again imply the existence of a shortest path from $v_0$ to the apex, hence also from $x_0$ to the apex, 
lying strictly to the left of the right boundary. The clockwise sequence of labels is thus necessarily $\ell-2\to\ell-1\to\ell-2\to\ell-3$. 
In particular, $w_0$ cannot be equal to $x_0$ since $w_0$ has a neighbor at distance $\ell-2$ which does not lie on the right boundary.
We conclude that there exists a path of two steps going from $v_0$ to a nearest neighbor $w_0$ at distance $\ell-1$ distinct from $x_0$ and then to a 
next-nearest neighbor $w_0'$ at distance $\ell-2$ distinct from $v_0$. Let us pick the \emph{leftmost such path of two steps}, i.e.\ going from 
$v_0$ to a nearest neighbor $x_1$ at distance $\ell-1$ distinct from $x_0$  
and then to a next-nearest neighbor $v_1$ at distance $\ell-2$ distinct from $v_0$. We may now
draw the leftmost shortest path $\mathcal{P}_1$ from $x_1$ to the apex, starting with the edge $x_1\to v_1$,  and call $v_1'$ the vertex at distance $\ell-3$
along $\mathcal{P}_1$. Considering the face immediately on the left of the edge of $\mathcal{P}_1$ from $v_1$ to $v_1'$,
and calling $w_1$ and $w_1'$ the two other incident vertices, again the four vertices $v_1,w_1,w_1',v_1'$ around
the face are necessarily distinct as otherwise, $\mathcal{P}_1$  would not be a leftmost shortest path, and, for the
same reasons as above, $w_1$ and $w_1'$
are necessarily at respective distances $\ell-1$ and $\ell-2$ from the apex. In particular, $w_1$ cannot be equal to $x_1$ as 
otherwise, $x_1$ would have a neighbor $w_1'$ at distance $\ell-2$ distinct from $v_0$ \footnote{The fact that $w_1'$ is itself distinct from $v_0$ is because otherwise, the edge $w_1'\to v_1'$ would 
lie strictly to the left of the right boundary and connect $v_0$ to a vertex at distance $\ell-3$, which is forbidden.}
and strictly to the left of $v_1$,
a contradiction.
To summarize, this proves the existence of a path of two steps
going from $v_1$ to a nearest neighbor at distance $\ell-1$ distinct from $x_1$ 
and then to a next-nearest neighbor at distance $\ell-2$ distinct from $v_1$. Again we pick the 
leftmost such path and call $x_2$ and $v_2$ the corresponding vertices (see figure~\ref{fig:dividing}), then drawn the leftmost shortest path $\mathcal{P}_2$ from $x_2$ to the apex (starting with the edge $x_2\to v_2$). Continuing this way, we build
by simple concatenation of the right-boundary edge from $x_0$ to $v_0$ and of all the elementary two-step paths a path connecting 
$x_0$ to $v_0$ to $x_1$ to $v_1$ to $x_2$ to $v_2$  
and so on, where all the $x_i$'s are at distance 
$\ell-1$ from the apex and all the $v_i$'s at distance $\ell-2$. 
\begin{figure}
\begin{center}
\includegraphics[width=11cm]{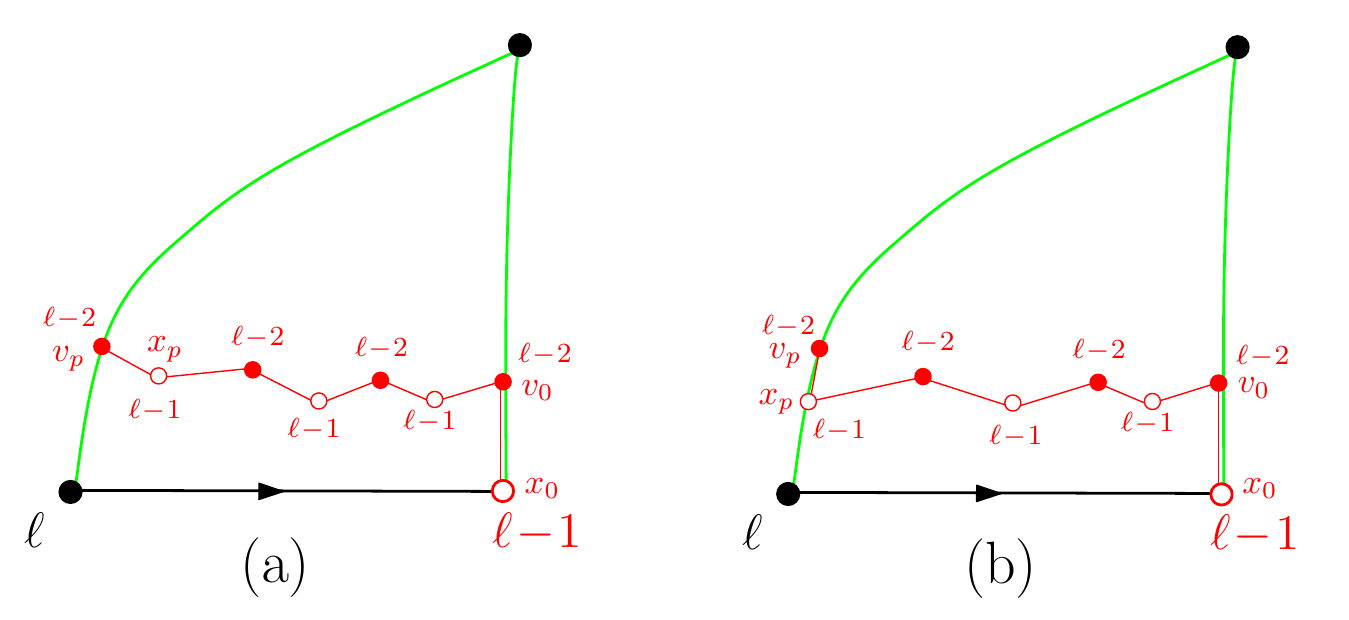}
\end{center}
\caption{The two possible ways for the dividing line to connect to the left boundary.}
\label{fig:twosituations}
\end{figure}
As explained below, this path \emph{cannot form a loop} so it defines an open
simple curve which necessarily reaches, after $p$ iterations of the process ($p\geq 1$), the vertex $v_p$ lying on the left boundary at distance
$\ell-2$ from the apex: this path defines our dividing line. Note that two situations may occur according to whether $x_p$ itself belongs to the left boundary or not (see figure~\ref{fig:twosituations}).

By construction, the dividing line is thus a simple open curve connecting $x_0$ to $v_p$ by
visiting alternatively vertices at distance $\ell-1$ and $\ell-2$ from the apex. The line therefore separates 
two domains in the slice, an upper part containing the apex and a lower part containing the root-vertex. Clearly, since a path from the apex to
the lower part must cross the dividing line, all the the vertices strictly inside the lower part are at distance at least $\ell-1$ from the apex.
\begin{figure}
\begin{center}
\includegraphics[width=6cm]{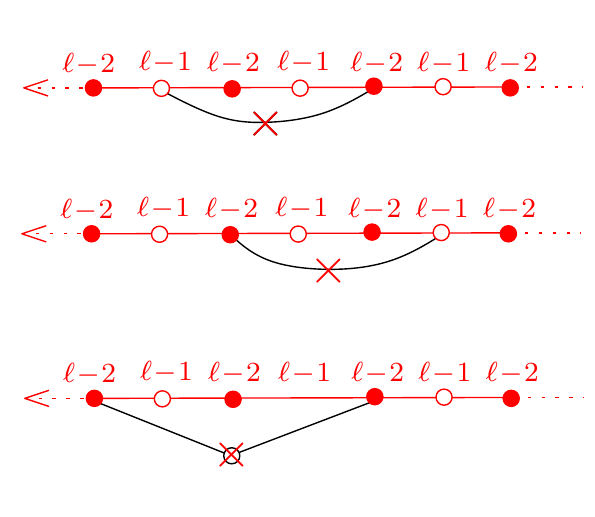}
\end{center}
\caption{A schematic picture of the properties of the dividing line emphasized in Property 1. Vertices at distance $\ell-2$ 
(resp. $\ell-1$) are colored in black (resp. in white).}
\label{fig:property1}
\end{figure}
By construction, the dividing line satisfies moreover the following property (illustrated in figure~\ref{fig:property1}):
\begin{proper}{} \hfill
\begin{itemize}
\item \emph{Two vertices of the dividing line cannot be linked by an edge lying strictly inside the lower part.}
\item \emph{Two vertices of the dividing line at distance $\ell-2$ cannot have a common neighbor strictly inside the lower part.}
\end{itemize}
\end{proper}
The first statement is clear as the existence of an edge linking two vertices of the dividing line and  inside the lower part 
would produce at some iteration of the dividing line construction an acceptable two-step path lying to the left of the chosen one.
As for the second statement, the common vertex would necessarily be at distance $\ell-1$ and again an acceptable two-step path 
would lie to the left of the chosen one. Note that, by contrast, pairs of vertices of the dividing line at distance $\ell-1$
may have a common neighbor strictly inside the lower part.

\begin{figure}
\begin{center}
\includegraphics[width=10cm]{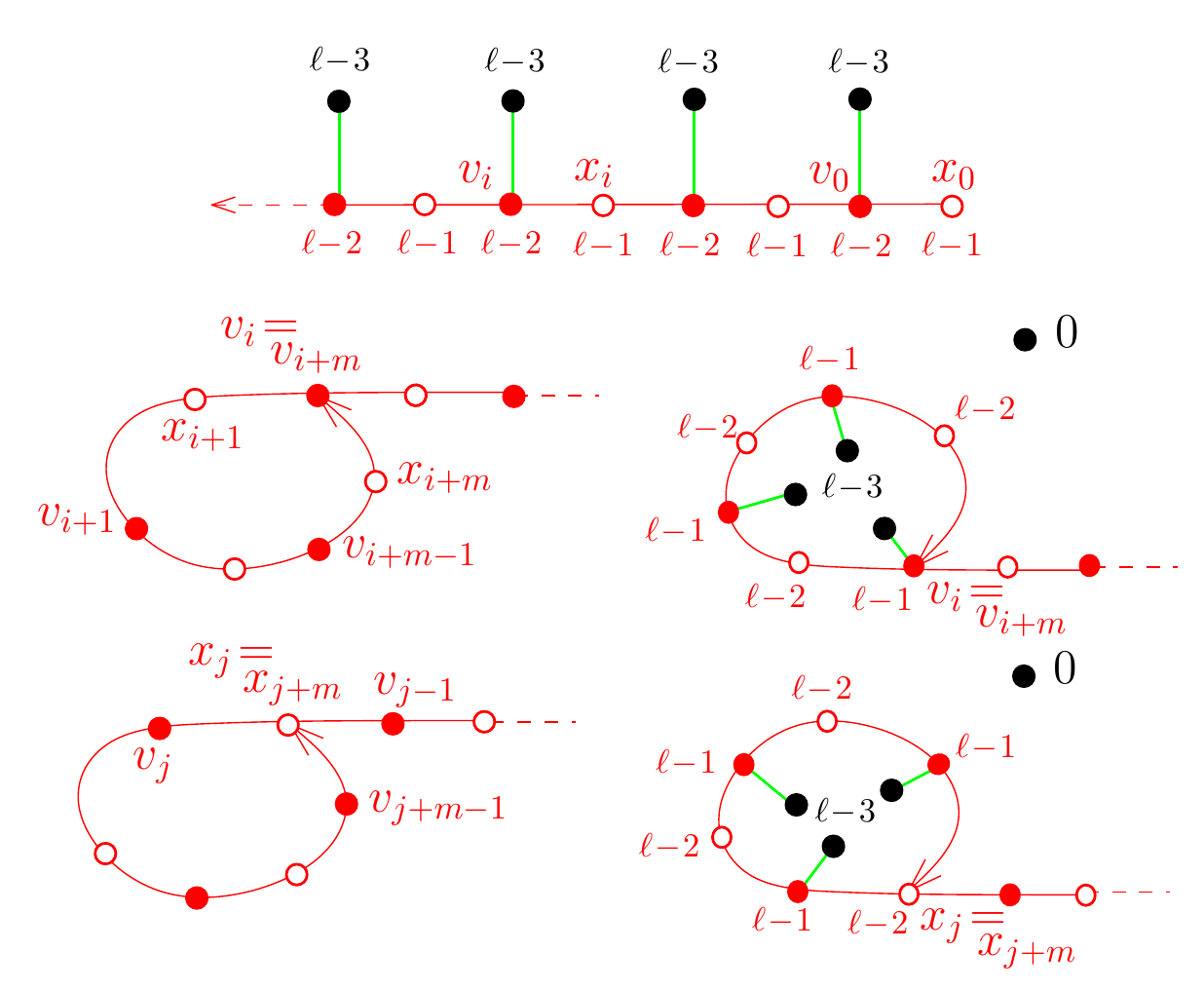}
\end{center}
\caption{A schematic picture of the dividing line (top, in red) made of its succession of vertices $x_i$ at distance $\ell-1$ and $v_i$ 
at distance $\ell-2$, with a vertex at distance $\ell-3$ attached to the right of each $v_i$. The lower part of the figure illustrates the contradictions 
which would occur if the dividing line were making a loop (see the arguments in the text).}
\label{fig:noloop}
\end{figure}
The fact that the concatenation of our two-step paths cannot form a loop may be understood via arguments
similar to those discussed in \cite{G15}. The proof is as follows: assume that the line forms a loop and 
consider the first vertex $v_i$, $i\geq 0$ (or respectively $x_j$, $j>0$) at which a double point arrises, i.e.\ $v_{i+m}=v_i$ 
for some $m>1$ (respectively $x_{j+m}=x_j$). Note that $m=1$ is not possible from our construction of the two-step paths.
Note also that the connection cannot occur at $x_0$ as this vertex has only one neighbor, $v_0$, at distance $\ell-2$.
If the connection occurs from the left (see figure~\ref{fig:noloop}), 
then the two-step path $v_i\to x_{i+m}\to v_{i+m-1}$ (respectively $v_{j-1}\to x_{j}\to v_{j+m-1}$) lies on the left of the chosen path $v_i\to x_{i+1}\to v_{i+1}$ (respectively $v_{j-1}\to x_{j}\to v_{j}$)
and should thus have been chosen instead of this latter path. This is a contradiction. If the connection occurs from the right, we use 
the property that, by construction, each vertex $v_n$ of the dividing line at distance $\ell-2$ from the apex has a neighbor on $\mathcal{P}_n$,
therefore \emph{on its right}, at distance $\ell-3$ from the apex (see figure~\ref{fig:noloop}). Then a loop closing from the right encloses at 
least one vertex 
at distance $\ell-3$ (for instance the neighbor of $v_{i+m}$, respectively of $v_{j+m-1}$) which is de facto surrounded by a frontier made of 
vertices at distance $\ell-2$ and $\ell-1$, a contraction.
We conclude that the dividing line cannot form a loop and necessarily ends on the left boundary.

As a final remark, we considered so far slices whose left-boundary length $\ell$ satisfies $\ell\geq 3$. When dealing with $T_k$, 
we also need to consider slices with $\ell=2$. For such slices, we define the dividing line as made of the single right-boundary edge linking
the endpoint $x_0$ of the root-edge to the apex $v_0$.   

\subsection{Decomposition of slices}
\label{sec:decomp}
As in \cite{G15}, the dividing line allows us to decompose slices enumerated by $T_k$ in a way which leads to a direct recursion
relating $T_k$ to $T_{k-1}$. As in previous section, the sequence of vertices along the dividing line will be denoted 
by $(x_0,v_0,x_1,v_1,\cdots,x_p,v_p)$
with $p\geq 0$, where the vertices $x_i$ (respectively $v_i$) are at distance $\ell-1$ (respectively $\ell-2$) from the 
apex, $\ell$ being the left-boundary length of the slice (and the constraint that $p=0$ if and only if $\ell=2$).
\begin{figure}
\begin{center}
\includegraphics[width=10cm]{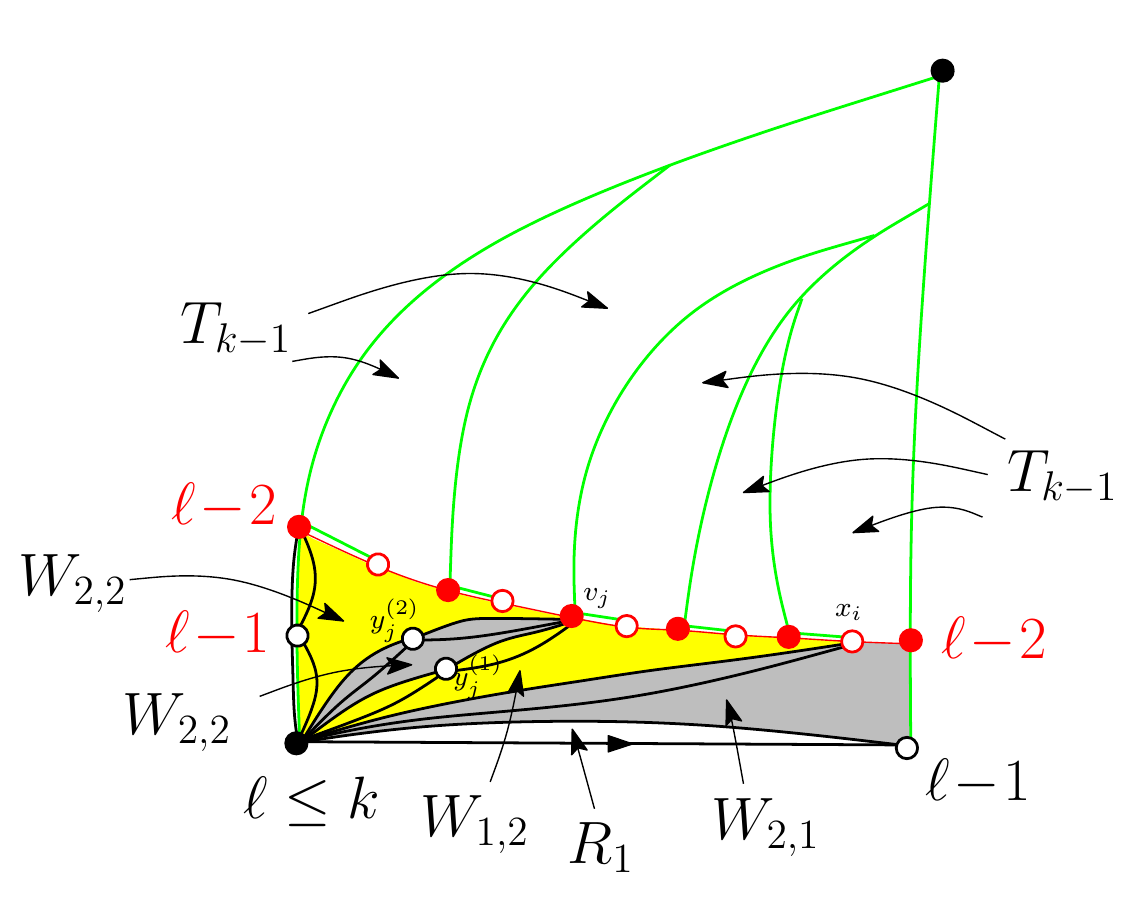}
\end{center}
\caption{The decomposition of a slice on both sides of the dividing line (see text). In the upper part, a slice enumerated 
by $T_{k-1}$ is created for each black-white edge of the dividing line (represented in red, supposedly oriented toward the left - here black and white vertices are represented 
by filled or empty red circles). In the lower part, a first bundle (enumerated by $R_1$) between the extremities of the root-edge is
completed by a sequence of blocks (here represented alternatively in grey and yellow). Each block has its right and left frontiers of lengths
$1$ or $2$ and is enumerated by $W_{1,1}$, $W_{1,2}$, $W_{2,1}$ or $W_{2,2}$ accordingly. In the present illustration, we 
have one $(2,1)$ block, one $(1,2)$ block and two $(2,2)$ blocks.}
\label{fig:Tknewrec}
\end{figure}
The decomposition is as follows: as already mentioned, the dividing line separates the slice 
into two domains, a lower part which contains the root-vertex and a complementary upper part (empty if and only if $\ell=2$). 
For $\ell\geq 3$, this upper part may be decomposed into slices by drawing the leftmost 
shortest path $\mathcal{P}_i$ from each vertex $x_i$ ($0\leq i\leq p$) to the apex (the path starting with  the edge of the dividing line linking 
$x_i$ to $v_i$). Note that $\mathcal{P}_0$ is the right boundary while $\mathcal{P}_p$ sticks to the left boundary from $v_p$ to the apex.
The paths $\mathcal{P}_i$ decompose the upper part in $p$ slices 
of left-boundary lengths between $2$ and $k-1$, hence enumerated by $T_{k-1}$, with a slice associated to each
step $v_{i-1}\to x_{i}$ ($1\leq i\leq p$) along the dividing line (the root-edge of this slice being the edge of the dividing line linking  $v_{i-1}$ to $x_{i}$, 
oriented from $x_{i}$ to $v_{i-1}$ -- se figure~\ref{fig:Tknewrec}).

More interesting is the decomposition of the lower part. We start by looking at the connections of the root-vertex to the dividing line:
the root vertex, at distance $\ell$ from the apex, is, in all generality, adjacent to a number of vertices $x_i$ of the dividing line.
These include $x_0$ plus possibly a number of other $x_i$'s, for instance $x_p$ in the situation (b) of figure \ref{fig:twosituations}.   
Note that these connections are in general achieved by a bundle (whose boundary if formed by the extremal 
edges performing the connection from the root-vertex to $x_i$). Now the root-vertex is also in all generality, connected
to a number of vertices $v_j$ of the dividing line by two-step paths whose intermediate vertex $y_j$ lies strictly 
inside the lower domain (and is at distance $\ell-1$ from the apex). These include for instance $v_p$ in the situation (a) of 
figure \ref{fig:twosituations}. The connection from the root-vertex to $y_j$ and from
$y_j$ to $v_j$ is achieved in general by a \emph{pair of bundles}. Moreover several intermediate vertices $y_j^{(1)},y_j^{(2)},\cdots, 
y_j^{(m)}$ may exist for the same $v_j$ (see figure ~\ref{fig:Tknewrec}). Now for each connection from the root-vertex to some $x_i$, we cut along the \emph{leftmost} edge
performing this connection and for each two-step-path connection from the root-vertex to some $v_j$ via some $y_j$, we cut along the
 \emph{leftmost} two-step path performing the connection. If several intermediate vertices $y_j^{(1)},y_j^{(2)},\cdots,
y_j^{(m)}$ exist, we make one cut for \emph{each occurrence} of such a vertex (see figure~\ref{fig:Tknewrec}). These cuts divide the lower part into
a sequence of connected domains whose left and right frontiers correspond to the performed cuts and have length
$1$ if the corresponding cut leads to  some $x_i$ or $2$ if the corresponding cut leads to some $v_j$. 
The domains may thus
be classified in four categories: $(1,1)$, $(1,2)$, $(2,1)$ or $(2,2)$ according to their right- and left-frontier length respectively
(for instance, the type $(1,2)$ corresponds to a right-frontier length $1$ and a left-frontier length $2$).
To be precise, the decomposition of the lower part may be characterized by some sequence $a_0,a_1,\cdots, a_n$, $n\geq 1$ (with $a_i,\in \{1,2\}$) corresponding to the successive 
encountered frontier lengths for the cut domains.
To each elementary step $a_i\to a_{i+1}$ of the sequence is attached a block of type $(a_i,a_{i+1})$.

The beginning of the sequence requires some special attention. Indeed, 
in the cutting, a first bundle from the root-vertex to $x_0$ is delimited, enumerated by $R_1$ (see figure ~\ref{fig:Tknewrec}),
Ignoring this first bundle, the effective right frontier of the first block is therefore a two-step path from the root-vertex to $y_0=x_0$,
then to $v_0$. We thus should start our sequence with $a_0=2$ (note that the root vertex may be also connected to $v_0$ by
two-step paths lying strictly inside the lower part, in which
case $a_1=2$ too). The sequence ends either with $a_n=2$ if the dividing line is in the situation (a) of figure \ref{fig:twosituations}  
or with $a_n=1$ if the dividing line is in the situation (b) of figure \ref{fig:twosituations}. 
\begin{figure}
\begin{center}
\includegraphics[width=11cm]{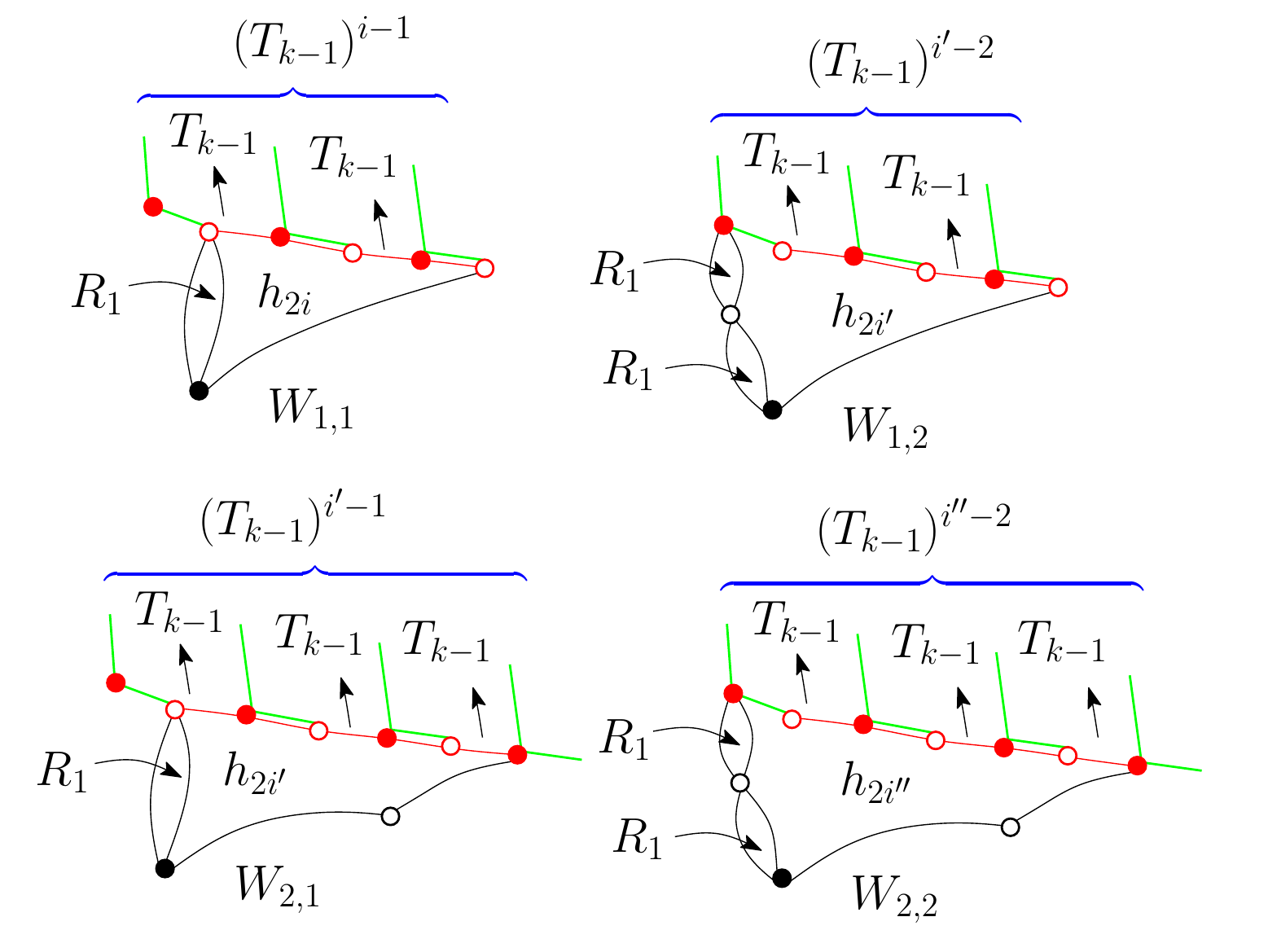}
\end{center}
\caption{Computation of the weights $W_{1,1}$, $W_{1,2}$, $W_{2,1}$ and $W_{2,2}$. For instance, the top-left situation 
corresponds to maps enumerated by $W_{1,1}$ for which 
we may extract a bundle (enumerated by $R_1$) along the left frontier. The remaining part, of length $2i$ (here $i=3$) is
enumerated by $h_{2i}$ and gives rise to $i-1$ black-white (from right to left) edges on the dividing line, hence $i-1$ factors $T_{k-1}$.
The other situations are represented with $i'=4$ (for $W_{1,2}$, and $W_{2,1}$) and $i''=5$ (for $W_{2,2}$).}
\label{fig:Wweights}
\end{figure}

Let us now discuss the weight that should be attached to each block of the sequence if we wish to compute $T_k$. 
Consider first a block of type $(1,1)$ (see figure~\ref{fig:Wweights}): it has an overall boundary of length $2i\geq 4$ made of its right
frontier (a single edge of length $1$), its left frontier (a single edge of length $1$ which is in general the leftmost edge of a bundle) and
a portion of length $2i-2$ of the dividing line which goes from some $x_j$ to $x_{j+i-1}$, hence contains $i-1$ edges of type  
$v_m\to x_{m+1}$, giving rise to $i-1$ slices in the upper part, hence producing a weight $T_{k-1}^{i-1}$. As for the block itself in the lower part,
we may decide to cut out the bundle to which belongs the left frontier of the block, giving rise to a weight $R_1$. The remaining
part (which has now as left frontier the rightmost edge of the bundle) is enumerated by some generating function $h_{2i}=h_{2i}(g)$ for particular quadrangulations with a boundary of length $2i$ satisfying special constraints which we will discuss below
(see Property 2). 
At this stage, let us just mention that we decide to choose as root-edge for these
quadrangulation the edge starting from the root-vertex counterclockwise around the domain. 

To summarize, the weight attached to a $(1,1)$ block is
\begin{equation*}
W_{1,1}=R_1\, \sum_{i\geq 2} h_{2i} \, T_{k-1}^{i-1}\ .
\end{equation*}
If we now consider a block of type $(1,2)$ of the lower domain, it has an overall boundary of length $2i\geq 4$ made of its right
frontier (a single edge of length $1$), its left frontier (a two-step path of length $2$ which is in general part of a pair of bundles)
a portion of length $2i-3$ of the dividing line which goes from some $x_j$ to $v_{j+i-2}$, hence contains $i-2$ edges of type  
$v_m\to x_{m+1}$, giving rise to $i-2$ slices in the upper domain, hence a weight $T_{k-1}^{i-2}$. As for the part in the lower domain,
we again decide to cut out the pair of bundles to which belongs the left frontier of the block, giving rise to a weight $R_1^2$. The remaining
part is again enumerated by $h_{2i}$ since the remaining quadrangulation of boundary length $i$ is precisely of the same type as above. 
The weight attached to a $(1,2)$ block is therefore
\begin{equation*}
W_{1,1}=R_1^2\, \sum_{i\geq 2} h_{2i} \, T_{k-1}^{i-2}\ .
\end{equation*}
Repeating the argument for $(2,1)$ and $(2,2)$ blocks, we find (see figure~\ref{fig:Wweights})
\begin{equation*}
W_{2,1}=W_{1,1}\, \qquad W_{2,2}=W_{1,2}\ ,
\end{equation*}
so that $W_{a,a'}$ actually depends only on the second index $a'$ (this is because both the number of bundles on the left side of the block
and the number of created slices in the upper part for a fixed $i$ depend only on $a'$ -- see figure~\ref{fig:Wweights}).
To get $T_k$, we must sum over all possible sequences $a_0,a_1,\cdots, a_n$. Since $a_0=2$ is fixed and all the other $a_i$'s are 
free (including $a_n$) and since the weights $W_{a,a'}$ depends only on $a'$, we immediately deduce the contribution
\begin{equation*}
R_1\, \sum_{a_0=2\atop a_1, \cdots, a_n \in \{1,2\}}\prod_{i=0}^{n-1}W_{a_i,a_{i+1}}=R_1\, \left(R_1\, \sum_{i\geq 2} h_{2i} \, T_{k-1}^{i-1}+R_1^2\, \sum_{i\geq 2} h_{2i} \, T_{k-1}^{i-2}\right)^n
\end{equation*}
for a sequence of length $n$, where we re-introduced the weight $R_1$ for the bundle from the root-vertex to $x_0$.
Recall that $n\geq 1$ since we have at least one block. Note also that slices with $\ell=2$ contribute to all values of $n$
via the $i=2$ term of the second sum in the parenthesis above\footnote{For $\ell=2$, we have
$a_i=2$ for all $i$ ($0\leq i\leq n$) and all the blocks have boundary length $4$, hence are enumerated by $R_1^2\, h_4$. 
In particular with have $T_2=R_1\sum_{n\geq 1}(R_1^2 h_4)^n=R_1^3\, h_4/(1-R_1^2\, h_4)$, in agreement with \eqref{eq:newrecur}
for $k=2$, since $T_1=0$ and $\Phi(0)=h_4$.}.

Summing over all $n\geq 1$, we arrive at the desired recursion relation
\begin{equation*}
T_k=\frac{R_1\, \left(R_1\, \sum_{i\geq 2} h_{2i} \, T_{k-1}^{i-1}+R_1^2\, \sum_{i\geq 2} h_{2i} \, T_{k-1}^{i-2}\right)}{1-
\left(R_1\, \sum_{i\geq 2} h_{2i} \, T_{k-1}^{i-1}+R_1^2\, \sum_{i\geq 2} h_{2i} \, T_{k-1}^{i-2}\right)}
\end{equation*}
or in short
\begin{equation}
T_k=\frac{R_1^2 (T_{k-1}+R_1)\, \Phi(T_{k-1})}{1-R_1(T_{k-1}+R_1)\, \Phi(T_{k-1})}\ , \qquad \Phi(T)\equiv \Phi(T,g)=\sum_{i\geq 2} h_{2i}(g)\, T^{i-2}\ .
\label{eq:newrecur}
\end{equation}
\begin{figure}
\begin{center}
\includegraphics[width=4cm]{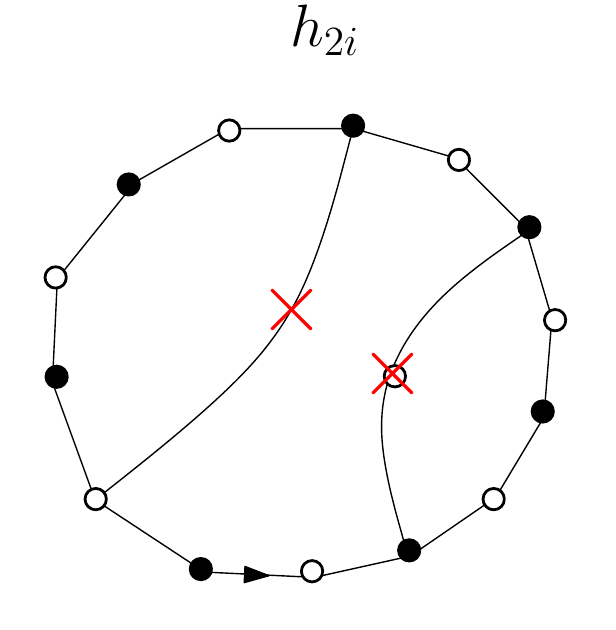}
\end{center}
\caption{A schematic representation of the two forbidden connections within maps enumerated by $h_{2i}$, as
dictated by Property 2.}
\label{fig:property2}
\end{figure}
To end the section, it remains to characterize the quadrangulations with a boundary of length $2i$ ($i\geq 2$) enumerated by $h_{2i}\equiv h_{2i}(g)$,
where $g$ is the weight per face.
By construction, the boundary of these quadrangulations is a simple curve and we may for convenience decide to color the boundary vertices alternatively in black and white, the root-vertex being black (thus the $v_i$'s at hand are also black and the $x_j$ and $y_j^{(\alpha)}$
are white). The maps enumerated by $h_{2i}$ are further characterized by the following property (see figure~\ref{fig:property2}):
\begin{proper} \hfill
\begin{itemize}
\item \emph{In the maps enumerated by $h_{2i}$, two vertices of the boundary cannot be linked by an edge lying strictly inside the map.}
\item \emph{In the maps enumerated by $h_{2i}$, two black vertices of the boundary cannot have a common (white) adjadent vertex strictly inside 
the map.}
\end{itemize}
\end{proper}
For pairs of vertices belonging to the dividing line, these properties are a direct consequence of Property 1. 
For pairs involving the other vertices (i.e.\ the root vertex or the intermediate vertices $y_j$), these properties are a direct
consequence of the block decomposition. As in \cite{G15}, it is remarkable that, while, in the decomposition, boundary vertices 
of the domains enumerated by $h_{2i}$ play different roles, the characterization of these domains via Property 2 turns out
to be symmetric for all boundary vertices.

\section{Simple quadrangulations}
\label{sec:simplemaps}
\subsection{From general to simple quadrangulations}
\label{sec:substitution}
As in \cite{G15}, we may slightly simplify our recursion by eliminating $R_1$ from our problem.
At the level of maps, it amounts to restrict our analysis to \emph{simple} quadrangulations (with a boundary),
i.e. quadrangulations \emph{without multiple edges}. We thus define simple analogs of $R_k$ and $T_k$,
namely the generating function $r_k\equiv r_k(G)$ of simple slices with left-boundary length $\ell$ in the range $1\leq \ell\leq k$
and $t_k\equiv t_k(G)$ for simple slices with $2\leq\ell\leq k$, with a weight $G$ per  inner face. Similarly, we define 
$\tilde{h}_{2i}\equiv \tilde{h}_{2i}(G)$ as the generating function, with a weight $G$ per inner face, 
of simple quadrangulations with a boundary of length $2i$ forming a simple curve (i.e.\ which does not cross itself),
and which satisfy Property 2. As it is well-known, we may pass from simple quadrangulations to general quadrangulations
by a substitution in the generating functions. Indeed, a general quadrangulation is obtained from a simple one by
replacing each edge of the simple quadrangulation by a bundle, as we defined it. Since the generating function for bundles
is $R_1$, the generating functions $R_k$, $T_k$, and $h_{2i}$ may in practice be obtained from $r_k$, $t_k$, and $\tilde{h}_{2i}$
by a substitution as follows: consider a quandrangulation with $F$ inner faces, $E$ inner edges and $2L$ edges 
on the boundary. We have the relation $4F=2E+2L$ 
so that
 \begin{equation*}
 E=2F-L\ .
 \end{equation*}  
In the case of maps enumerated by $R_k$ (respectively $T_k$), we must, starting from maps enumerated by $r_k$
(respectively $t_k$), put a weight $R_1$ to each 
inner edge as well as to each of the $L=\ell$ edges of the left boundary and finally to the root edge. No weight $R_1$ is assigned to the edges
of the right boundary as bundles cannot be present there since the right boundary is the unique shortest path between its 
extremities. We must thus assign a global weight $R_1^{E+L+1}=R_1 \times R_1^{2F}$ (note that, written this way, the weight is independent of $\ell$). 
In other words, we must assign a weight 
$R_1^2$ per face and a global factor $R_1$, which yields the relations
\begin{equation}
R_k(g)=R_1\, r_k(G)\ , \qquad T_k(g)=R_1\, t_k(G)\ ,
\label{eq:rRtT}
\end{equation}
with the correspondence
\begin{equation}
G=g\, R_1^2\ .
\label{eq:gG}
\end{equation}
Note that the relation $T_k=R_k-R_1$ translates into
\begin{equation*}
t_k=r_k-1
\end{equation*} 
consistent with the fact that there is a unique simple slice with $\ell=1$, the root-edge-map, hence $r_1=1$. As for $h_{2i}$, it is obtained 
by assigning a weight $R_1$ to all the inner edges of the maps enumerated by $\tilde{h}_{2i}$. Again, because of Property 2, 
there cannot be multiple edges on the boundary. We must thus assign a global weight $R_1^{E}=R_1^{2F-L}
=R_1^{2F-i}$ since $L=i$ in this case.  In other word, we must again assign an extra weight 
$R_1^2$ per face and now a global factor $R_1^{-i}$, which yields
\begin{equation}
h_{2i}(g)=R_1^{-i}\, \tilde{h}_{2i}(G)
\label{eq:hhtilde}
\end{equation}
with the correspondence \eqref{eq:gG}.
Introducing the quantity
\begin{equation}
\tilde{\Phi}(t)\equiv \tilde{\Phi}(t,G)=\sum_{i\geq 2} \tilde{h}_{2i}(G)\, t^{i-2}\ ,
\label{eq:defPhitilde}
\end{equation}
we deduce from \eqref{eq:hhtilde}
\begin{equation*}
\Phi(T)=R_1^{-2}\, \tilde{\Phi}(t)\ , \qquad T=R_1\, t
\end{equation*}
with the implicit correspondence \eqref{eq:gG}.

Finally, the recursion \eqref{eq:newrecur} translates into the simpler relation
\begin{equation}
t_k=\frac{(t_{k-1}+1)\, \tilde{\Phi}(t_{k-1})}{1-(t_{k-1}+1)\, \tilde{\Phi}(t_{k-1})}
\label{eq:newrecurbis}
\end{equation}  
(with $t_1=0$) where, as promised, $R_1$ is no longer present. As in \cite{G15}, we note that there is no straightforward analog
of the relation \eqref{eq:GkRk} for simple quadrangulations. This is because closing a slice into a planar quadrangulations 
by identifying its right and left boundaries as in figure~\ref{fig:twopoint} may in general create multiple edges. The recourse to simple slices in this paper should therefore simply be viewed as a non-essential but convenient way to slightly simplify our recursion by temporarily removing 
the $R_1$ factors.

\subsection{An equation for $\mathbf{\tilde{\Phi}(t)}$}
\label{sec:eqforPhi}
In order to solve \eqref{eq:newrecurbis}, and eventually \eqref{eq:newrecur}, we need some more explicit
expression for $\tilde{\Phi}(t)$. Such expression may be obtained by first noting that $\tilde{\Phi}(t)$
is fully determined by the following equation:
\begin{equation}
\tilde{\Phi}(t)=G+\frac{G}{t}\left\{\frac{(t+1)\tilde{\Phi}(t)}{1-(t+1)\tilde{\Phi}(t)}-\frac{\tilde{h}_4}{1-\tilde{h}_4}\right\}
\label{eq:eqforPhitilde}
\end{equation}  
which may equivalently be written as
\begin{equation*}
t(1+t)\tilde{\Phi}^2(t)+\left(G(1+t)(1-t+g_4)-t\right)\tilde{\Phi}(t)+G(t-g_4)=0\ , \qquad g_4\equiv \frac{\tilde{h}_4}{1-\tilde{h}_4}\ .
\end{equation*}
\begin{figure}
\begin{center}
\includegraphics[width=8cm]{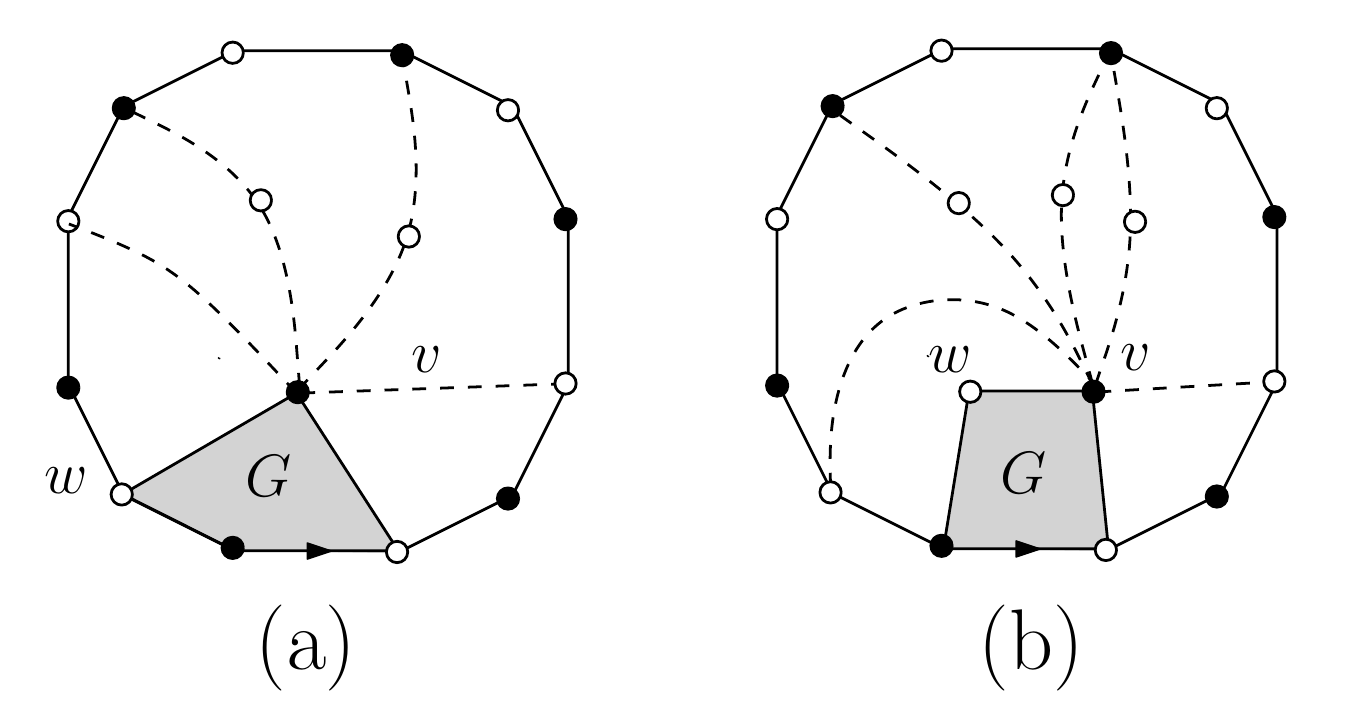}
\end{center}
\caption{The decomposition of a map enumerated by $\tilde{\Phi}(t)$ (and not reduced to a single face). The map is decomposed by
removing the face (in gray -- weight $G$) on the left of the root-edge and cutting along all the connections of the black vertex $v$ (incident to the gray
face and strictly inside the map) to boundary vertices
 by either simple edges or by two-step paths. The white vertex $w$ (incident to the gray face and different from the root-edge 
 extremities) may lie on the boundary (case (a)) or not (case (b)). In this latter case, it cannot be connected to other black boundary
 vertices than the root-vertex.}
\label{fig:Phidecomp}
\end{figure}
\begin{figure}
\begin{center}
\includegraphics[width=4.5cm]{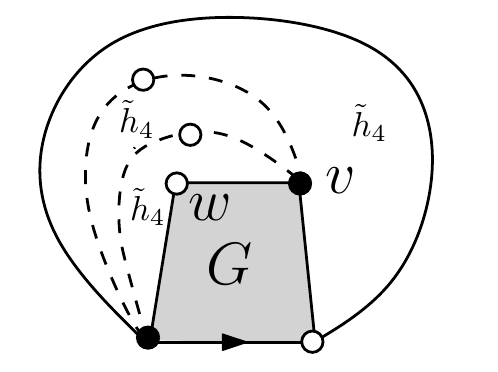}
\end{center}
\caption{A schematic picture of those maps which have the structure of figure~\ref{fig:Phidecomp} but which do not contribute to $\tilde{\Phi}(t)$ since their boundary length is $2$
(hence less than $4$).}
\label{fig:excluded}
\end{figure}

Let us first prove \eqref{eq:eqforPhitilde} and then show how to get $\tilde{\Phi}(t)$ out of it. From the definition \eqref{eq:defPhitilde}, 
$t^2\, \tilde{\Phi}(t)$
enumerates simple quandrangulations with a boundary of \emph{arbitrary length $2i$} ($i\geq 2$) forming a simple curve
and satisfying Property 2, with a weight $G$ per inner face and a weight $\sqrt{t}$ per boundary edge.
The quadrangulation may be reduced to a single face (with a boundary of length $4$), leading to a first contribution $t^2\,  G$ to $t^2\, \tilde{\Phi}(t)$,
hence (after dividing by $t^2$) to the first term in \eqref{eq:eqforPhitilde}.
In all the other cases, we may look at the face immediately on the left of the root-edge and call $v$ and $w$ 
its black and white incident vertices other than the extremities of the root-edge (recall that we decided 
for convenience to bi-color the map in black and white, the root-vertex being black). The (black) vertex $v$ 
\emph{cannot lie on the boundary} as otherwise, $w$ (which could not lie in this case on the boundary because
of the first requirement of Property 2) would be a common neighbor to $v$ and to the root-vertex, thus violating
the second requirement of Property 2. So $v$ lies strictly inside the map.
As for the (white) vertex $w$, it may lie on the boundary (see figure~\ref{fig:Phidecomp}-case (a)) or not
(case (b)). In this latter case, $w$ \emph{cannot be connected to a vertex of the boundary} other than the root-vertex
as otherwise,   the second requirement of Property 2 would again be violated.
On the other hand, the vertex $v$ is connected by  simple edges to a number of white vertices 
of the boundary (including the endpoint of the root-edge as well as $w$ in case (a)), and may be connected 
by two-step paths to black vertices of the boundary (including the root-vertex in case (b)). Let us draw all these
connections and cut the map along them. After cutting, the face to the left of the root edge gets disconnected
and contributes a weight $t\, G$ to $t^2\, \tilde{\Phi}(t)$ in case (a) and a weight $\sqrt{t}\, G$ in case (b).
The rest of the map forms a sequence of blocks. As we did before in the slice decomposition, we may consider 
the sequence $a_0,a_1, \cdots, a_{n}$ ($n\geq 1$) of the lengths $1$ or $2$ of the (counterclockwise) successive connections 
of $v$ to boundary vertices (with $1$ for a simple edge connection to a white vertex and $2$ for a two-step-path
connection to a black vertex). We have $a_0=1$ since the first connection is from $v$ to the white endpoint of the root-edge
while $a_{n}=1$ in case (a) and $a_{n}=2$ in case (b). Now the $m$-th block has a boundary of total arbitrary length $2j$ for
some $j\geq 2$, with $2j-a_{m-1}-a_{m}$ edges on the original boundary of the map. It must thus be given a 
weight $w(a_{m-1},a_m)$ with 
\begin{equation*}
\hskip -.8cm w(a,a')=\sum_{j\geq 2} \tilde{h}_{2j} (\sqrt{t})^{2j-a-a'} = (\sqrt{t})^{4-a-a'}\, \tilde{\Phi}(t) 
\end{equation*}  
Considering both cases (a) and (b), the contribution to $t^2\, \tilde{\Phi}(t)$ of sequences of $n$ blocks (together 
with that of the face on the left of the root-edge) is therefore
\begin{equation*}
\begin{split}
&\hskip -1.cm \left(\tilde{\Phi}(t)\right)^n \left( \sum_{a_0=1, a_n=1 \atop a_2,\cdots, a_{n-1}\in \{1,2\}}\!\!\!\!\!\!\!\!\! t\, G\, \prod_{m=1}^{n} (\sqrt{t})^{4-a_{m-1}-a_{m}}
\quad + \sum_{a_0=1, a_n=2 \atop a_2,\cdots, a_{n-1}\in \{1,2\}}\!\!\!\!\!\!\!\!\! \sqrt{t} \, G \, \prod_{m=1}^{n}  (\sqrt{t})^{4-a_{m-1}-a_{m}}\right)\\
&\hskip -1.cm = t\, G\, \left(\tilde{\Phi}(t)\right)^n \left( \sum_{a_n=1 \atop a_2,\cdots, a_{n-1}\in \{1,2\}}  \prod_{m=1}^{n} t^{2-a_m}
\quad + \sum_{a_n=2 \atop a_2,\cdots, a_{n-1}\in \{1,2\}} \prod_{m=1}^{n}  t^{2-a_m}\right)\\
&\hskip -1.cm = t\, G\, \left(\tilde{\Phi}(t)\right)^n  \sum_{a_1,a_2,\cdots, a_{n}\in \{1,2\}}  \prod_{m=1}^{n} t^{2-a_m}\\
&\hskip -1.cm = t\, G\, \left(\tilde{\Phi}(t)\right)^n (t+1)^n\ .\\
\end{split}
\end{equation*}
To go from the first to the second line, we simply use the fact that $a_0=a_n$ for sequences of the first sum and
$a_0=a_n-1$ for sequences of the second sum. 
Summing over $n\geq 1$ yields the contribution 
\begin{equation}
t\, G\, \frac{(t+1)\tilde{\Phi}(t)}{1-(t+1)\tilde{\Phi}(t)}
\label{eq:secondcontrib}
\end{equation}
to $t^2\, \tilde{\Phi}(t)$, hence (after dividing by $t^2$) to the second term in \eqref{eq:eqforPhitilde}. 
So far in our decomposition, we did not enforce the condition that, in $\tilde{\Phi}(t)$, the length $2i$ of the
maps must satisfy $i\geq 2$. In our block sequences, there is a situation where this length happens to be $2$ (i.e.\ $i=1$)
(see figure~\ref{fig:excluded}):
it corresponds to a situation of case (b) with a first block of boundary length $4$ whose boundary vertices
are $v$, $w$ and the two extremities of the root-edge (and with exactly $1$ edge
on the original boundary of the whole map), completed by arbitrarily many blocks of size $4$ whose boundaries are made of 
two-step-paths from $v$ to the root-vertex (these blocks do not contribute to the original boundary length).  
These maps, made of $n\geq 1$ blocks enumerated by $\tilde{h}_4$ contribute
\begin{equation*}
t\, G\, \sum_{n\geq 1} (\tilde{h}_4)^n = t\, G\, \frac{\tilde{h}_4}{1-\tilde{h}_4}
\end{equation*}
to \eqref{eq:secondcontrib} and must be subtracted to properly recover $t^2\, \tilde{\Phi}(t)$.
This explains (after dividing by $t^2$) the third term in \eqref{eq:eqforPhitilde}.

\subsection{An expression for $\mathbf{\tilde{\Phi}(t)}$}
\label{sec:kernel}
We shall now extract from \eqref{eq:eqforPhitilde} a tractable expression for $\tilde{\Phi}(t)$. The first 
step consists in getting from the equation an expression for $\tilde{h}_4=\tilde{h}_4(G)$ as a function of
 the face weight $G$. Here we use the following standard trick: from \eqref{eq:eqforPhitilde},
we may write 
\begin{equation}
\tilde{h}_4=\frac{t(t+1)\tilde{\Phi}^2(t)-(G(t+1)(t-1)+t)\tilde{\Phi}(t)+G\, t}{t(t+1)\tilde{\Phi}^2(t)-(G(t+1)t+t)\tilde{\Phi}(t)+G\, (t+1)}
\label{eq:h4}
\end{equation}
which, upon differentiating with respect to $t$ (recall that $\tilde{h}_4$ does not depend on $t$), yields
\begin{equation*}
\begin{split}
0& =\frac{d\tilde{h}_4}{dt}\\ & \Rightarrow 
0 =\tilde{\Phi}'(t)\left\{
t\left((t+1)^2 \tilde{\Phi}^2(t)-2(t+1)\tilde{\Phi}(t)+1-G\right)-G
\right\}\\
&\ \ 
+\left\{ 
\tilde{\Phi}(t) \left((t+1)^2 \tilde{\Phi}^2(t)-2(t+1)\tilde{\Phi}(t)+1-G\right)-G\left(1-(t+1)\tilde{\Phi}(t)\right)^2
\right\}\ .\\
\end{split}
\end{equation*}
This equation is satisfied in particular if we let $t$ vary on a line $t=t(G)$ where each of the two terms between brackets
in the above expression vanishes.
Canceling these two terms and solving for $G$ and $\tilde{\Phi}$ yields the following two possible solutions
\begin{equation*}
\begin{split}
& G=\frac{t(G)}{(1+t(G))^3} \qquad {\rm and} \qquad \tilde{\Phi}(t(G))=\frac{t(G)}{(1+t(G))^2}\\
& {\rm or}\\
& G=\frac{1}{t(G)(1+t(G))} \qquad {\rm and} \qquad \tilde{\Phi}(t(G))=\frac{1}{t(G)}\\
\end{split}
\end{equation*}
which select two lines $t=t(G)$ where we know the value of $\tilde{\Phi}$. 
Plugging these values in \eqref{eq:h4} yields
\begin{equation*}
\tilde{h}_4=\frac{t(G)(1-t(G))}{1+t(G)-t(G)^2}
\end{equation*}
for the first choice and $\tilde{h}_4=1+G$ for the second choice. This latter result is clearly not satisfactory
(recall that $\tilde{h}_4$ enumerates simple quandrangulations with a boundary of length $4$ satisfying Property 2,
with a weight $G$ per face)
so we are left with the first choice. In other words we deduce from our particular solution the parametric expression 
for $\tilde{h}_4$:
\begin{equation}
\tilde{h}_4= \frac{C(1-C)}{1+C-C^2}\ ,\qquad {\rm where}\ G=\frac{C}{(1+C)^3} 
\label{eq:param}
\end{equation}
(here $C=t(G)$ should be viewed as a simple parametrization of $G$) which implicitly determines $\tilde{h}_4$ as a function of $G$.

Inverting the relation between $G$ and $C$, we get
\begin{equation}
C=\sum_{n\geq 1} \frac{1}{2n+1}{3n \choose n} G^n\ .
\label{eq:CG}
\end{equation}
\begin{figure}
\begin{center}
\includegraphics[width=12cm]{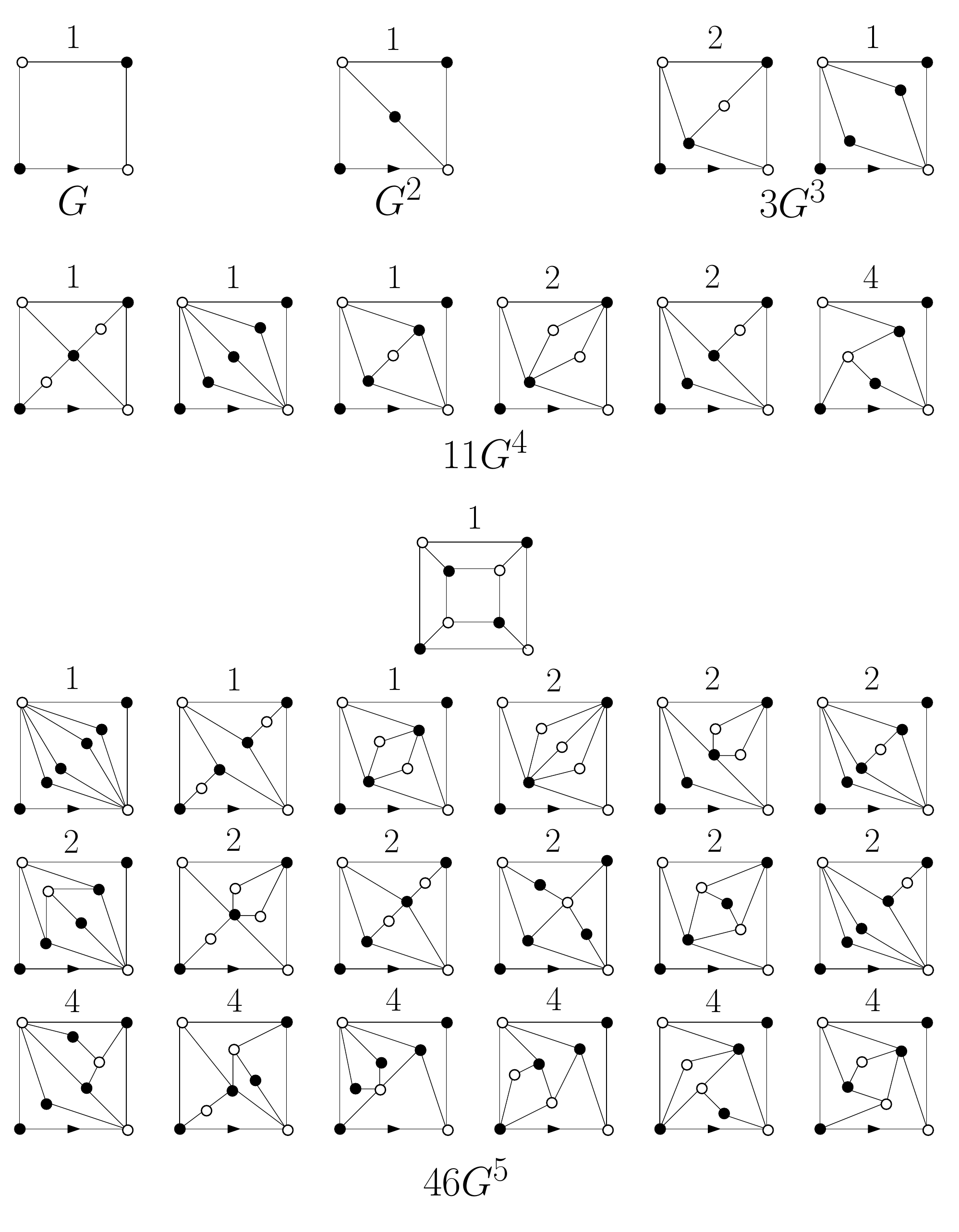}
\end{center}
\caption{A direct inspection of the maps enumerated by $\tilde{h}_4$ and with up to $5$ inner faces.
To obtain all the acceptable maps, each of the maps presented here must also be reflected along its two diagonals (and re-rooted at its bottom edge).
According to the symmetry of the map at hand, the number of \emph{distinct} maps obtained by these reflections is
$1$, $2$ or $4$ as indicated. Summing these numbers yields respectively $1,1,3,11$ and $46$  distinct maps with $1,2,3,4$ and $5$ inner faces,
explaining the first five terms of \eqref{eq:firstterms}.}
\label{fig:firstterms}
\end{figure}
Note that we implicitly assume that $0\leq G\leq 4/27$ for $C$ to be well-defined as above, which in turns
implies that $C$ lies in the range $0\leq C\leq 1/2$. This is consistent with the fact that the number of simple quandrangulations
grows exponentially with the number $F$ of faces as $(27/4)^F$ \cite{Brown1965}.
Now a Lagrange inversion yields the following formula 
\begin{equation*}
\hskip -1.2cm [G^p]\tilde{h}_4=\sum_{n=0}^{p-1}\left(\frac{(-1)^n\omega^{n+2}-\omega^{-n-2}}{\sqrt{5}}\right)\frac{(n+1)(3p)!}{p(p-1-n)!(2p+1+n)!}\ , \qquad \omega=\frac{1+\sqrt{5}}{2}\  
\end{equation*}
(involving Fibonacci numbers) for the number of simple quandrangulations with a boundary of length $4$, with $p$ inner faces,  satisfying Property 2. 
The first terms of this expansion read
\begin{equation}
\tilde{h}_4=G+G^2+3G^3+11G^4+46G^5+209G^6+\cdots
\label{eq:firstterms}
\end{equation}
with coefficients which may easily be verified by direct inspections of the maps at hand. (see figure \ref{fig:firstterms}).  
This sequence appears in \cite{Kuba} in the context of the enumeration of naturally embedded ternary trees. This should
not come as a surprise since bijections exist between such trees and simple quadrangulations \cite{JaSc98,BG10}.
\vskip .2cm
With the parametrization \eqref{eq:param}, equation \eqref{eq:eqforPhitilde} may be rewritten as a quadratic equation
\begin{equation}
\begin{split}
t(t+1)(1+C)^3\tilde{\Phi}^2(t)+& \left\{C(1\!+\!C\!-\!C^2)- t(1\!+\!3C\!+\!2C^2\!+\!2 C^3)-t^2\, C\right\}\tilde{\Phi}(t) \\
& \qquad \qquad \qquad \qquad \qquad \qquad \qquad +C(t-C+C^2)=0 \\
\end{split}
\label{eq:PhiC}
\end{equation}
whose discriminant reads
\begin{equation*}
\Delta=(C-t)^2 \left\{(1+C^2-C(t-1))^2- 4 C^2(1+C)(t+1)\right\}\ .
\end{equation*}
A look at the second factor suggests introducing the quantity\footnote{Here we use a trick in all points similar to that
used by Tutte in \cite{TutteCPT}.} $Y(t)$, solution of 
\begin{equation}
Y^2(t)+(1+C^2-C(t-1))\, Y(t)+C^2(1+C)(t+1)=0
\label{eq:defY}
\end{equation}
whose two solutions $Y_\pm(t)$ are related by the following involution (obtained by eliminating $t$ between the equation \eqref{eq:defY} for
$Y(t)=Y_+(t)$ and the same equation for $Y(t)=Y_-(t)$):
\begin{equation}
Y_\pm(t)=\frac{C(1+C)((1+C)^2+Y_\mp(t))}{Y_\mp(t)-C(1+C)}\ .
\label{eq:invol}
\end{equation}
For both determinations, we have the relation (directly read off \eqref{eq:defY}):
\begin{equation}
t=\frac{(1+C+Y(t))(C^2+Y(t))}{C(Y(t)-C-C^2)}\ .
\label{eq:tY}
\end{equation}
Plugging this value in \eqref{eq:PhiC} allows to rewrite the equation for $\tilde{\Phi}(t)$ as
\begin{equation*}
\begin{split}
\hskip -.5cm & \left(\tilde{\Phi}(t)-\frac{C(1+C+C^2+Y(t))}{(1+C+Y(t))((1+C)^2+Y(t))}\right)\\
&\qquad\qquad \times \left(\tilde{\Phi}(t)-
\frac{C(Y(t)-C(1+C))(Y(t)+C^2(1+C))}{(1+C)^3Y(t)(C^2+Y(t))}\right)=0\ .\\
\end{split}
\end{equation*}
This gives a priori two possible expressions for $\tilde{\Phi}(t)$ as a function of $Y(t)$ but
it is easily seen that the two formulas get interchanged by the involution \eqref{eq:invol}. We may therefore
decide to choose the expression coming from the first factor, namely
\begin{equation}
\tilde{\Phi}(t)=\frac{C(1+C-C^2+Y(t))}{(1+C+Y(t))((1+C)^2+Y(t))}
\label{eq:PhiY}
\end{equation}
provided we pick the correct determination of $Y(t)$. This determination is fixed by the small $t$ behavior
$\tilde{\Phi}(t)=\tilde{h}_4+O(t)$ with the formula \eqref{eq:param} for $\tilde{h}_4$.
This selects the determination (recall that $0\leq C\leq 1/2$) 
\begin{equation}
Y(t)=\frac{1}{2}\left\{C(t-1)-1-C^2+\sqrt{(1+C^2-C(t-1))^2- 4 C^2(1+C)(t+1)}\right\}
\label{eq:valY}
\end{equation}
(indeed, we then have $Y(0)=-C^2$ and we recover for $\tilde{\Phi}(0)$ the correct expression 
\eqref{eq:param} for $\tilde{h}_4$ while the other determination would yield $Y(0)=-C-1$ in which case 
$\tilde{\Phi}(t)$ would diverge for $t\to 0$).

Finally, from \eqref{eq:tY}, the expression \eqref{eq:PhiY} may be simplified into
\begin{equation}
\tilde{\Phi}(t)=\frac{C^2}{Y(t)(1+C+Y(t))}+\frac{1}{t+1}\ .
\label{eq:PhiYt}
\end{equation}
Equations \eqref{eq:PhiYt}, \eqref{eq:valY} and \eqref{eq:CG}, which implicitly fix $\tilde{\Phi}(t)=\tilde{\Phi}(t,G)$
in terms of $G$, 
are very reminiscent of similar expressions for triangulations (see \cite{TutteCPT,G15})

\section{Solution of the recursion}
\label{sec:final}
\subsection{Solution for simple quadrangulations}
\label{sec:simplesol}
We are now ready to solve \eqref{eq:newrecurbis}. Again, as in \cite{G15}, the idea is to rewrite this recursion in terms 
of the variable $Y(t)$. More precisely, let us define 
\begin{equation*}
Y_k\equiv Y(t_k)\ , \qquad \tilde{\Phi}_k\equiv\tilde{\Phi}(t_k)\ .
\end{equation*}
From \eqref{eq:tY}, we have the relations
\begin{equation}
t_k=\frac{(1+C+Y_k)(C^2+Y_k)}{C(Y_k-C-C^2)}\ , \qquad t_{k-1}=\frac{(1+C+Y_{k-1})(C^2+Y_{k-1})}{C(Y_{k-1}-C-C^2)}\ .
\label{eq:tkYk}
\end{equation}
Our recursion \eqref{eq:newrecurbis} reads
\begin{equation}
t_k=\frac{(t_{k-1}+1)\, \tilde{\Phi}_{k-1}}{1-(t_{k-1}+1)\, \tilde{\Phi}_{k-1}}
\label{eq:recurPhik}
\end{equation}
with, from \eqref{eq:PhiYt},
\begin{equation*}
\tilde{\Phi}_{k-1}=\frac{C^2}{Y_{k-1}(1+C+Y_{k-1})}+\frac{1}{t_{k-1}+1}\ .
\end{equation*}
Inserting this expression in \eqref{eq:recurPhik} and plugging the value \eqref{eq:tkYk} of $t_{k-1}$,
we obtain $t_k$ as a function of $Y_{k-1}$, namely:
\begin{equation*}
t_k=\frac{Y_{k-1}(C^2-C-1-Y_{k-1})}{C((1+C)^2+Y_{k-1})}\ .
\end{equation*} 
Equating this formula with that of \eqref{eq:tkYk} for $t_k$, we obtain the following relation
between $Y_k$ and $Y_{k-1}$:
\begin{equation*}
(1+C+Y_{k-1}+Y_k)\left(C^2(1+C)^2-C(1+C)Y_{k-1}+(1+C)^2 Y_k+Y_{k-1}Y_k\right)=0\ .
\end{equation*}
To choose which factor to cancel, we note that, for $G\to 0$, $t_k$ and $t_{k-1}$ tend to $0$ and thus $Y_k$ and $Y_{k-1}$
tend to $-C^2\to 0$
so the first factor does not vanish. The correct choice is therefore to cancel the second factor
and we arrive at the remarkably simple recursion relation for $Y_k$:
\begin{equation}
Y_k=\frac{C(1+C)\, Y_{k-1}-C^2(1+C)^2}{Y_{k-1}+(1+C)^2}\ .
\label{eq:recurYk}
\end{equation}
As recalled in \cite{G15}, solving such a recursion relation is a standard exercise. To solve more generally
the equation 
\begin{equation*}
Y_k=f(Y_{k-1})\ , \qquad f(Y)\equiv \frac{a\, Y+b}{c\, Y+d}\ ,
\end{equation*}
we introduce the two fixed points $\alpha$ and $\beta$ of the function $f$ (i.e.\ the two solutions of $f(Y)=Y$).
Then the quantity 
\begin{equation*}
W_k=\frac{Y_k-\alpha}{Y_k-\beta}
\end{equation*}
satisfies  $W_k=x \, W_{k-1}$, hence
\begin{equation*}
W_k=x^{k-1}\, W_1\ , \qquad x \equiv \frac{c\, \beta+d}{c\, \alpha+d}\ .
\end{equation*}
The desired $Y_k$ is recovered  
via $Y_k=(\alpha-\beta\, W_k)/(1-W_k)$ (note that $\alpha$ and $\beta$ are supposed to be distinct, as
will be verified a posteriori). 

In our case, we may take
\begin{equation*}
a=C(1+C)\ , \qquad b= - C^2(1+C)^2\ , \qquad c=1\ , \qquad d=(1+C)^2
\end{equation*}
so that 
\begin{equation*}
\begin{split}
\alpha & = \frac{1}{2}(1+C)\left(\sqrt{1-4 C^2}-1\right)\\
\beta & = \frac{1}{2}(1+C)\left(-\sqrt{1-4 C^2}-1\right)\\
x & = \frac{1-\sqrt{1-4 C^2}}{2C}\\
\end{split}
\end{equation*}
(in particular, $\alpha\neq\beta$ for $C\neq 1/2$, i.e.\ $x\neq 1$). 
The last formula is inverted into
\begin{equation*}
C=\frac{x}{1+x^2}
\end{equation*}
and the first two may then be rewritten as
\begin{equation*}
\begin{split}
\alpha & = -x^2\ \frac{1+x+x^2}{(1+x^2)^2}\\
\beta & = - \frac{1+x+x^2}{(1+x^2)^2} \ .\\
\end{split}
\end{equation*}
The initial condition reads
\begin{equation*}
t_1=0\ \Rightarrow \ Y_1=-C^2 \ \Rightarrow  \ W_1=\frac{-C^2-\alpha}{-C^2-\beta}=x^3
\end{equation*}
so that 
\begin{equation*}
W_{k}=x^{k-1}\, W_1= x^{k+2}
\end{equation*}
and 
\begin{equation*}
Y_k=- x^2 \ \frac{1+x+x^2}{(1+x^2)^2}\ \frac{1-x^k}{1-x^{k+2}}\ .
\end{equation*}
Plugging this expression in \eqref{eq:tkYk} gives
\begin{equation*}
t_k=\frac{x}{1+x^2}\ \frac{(1-x^{k-1})(1-x^{k+4})}{(1-x^{k+1})(1-x^{k+2})}
\end{equation*}
and eventually
\begin{equation*}
r_k=t_k+1=\frac{1+x+x^2}{1+x^2}\ \frac{(1-x^{k})(1-x^{k+3})}{(1-x^{k+1})(1-x^{k+2})}\ .
\end{equation*}
From the connection \eqref{eq:param} between $G$ and $C$ and that just above between $C$ and $x$, we deduce
the relation between $G$ and $x$:
\begin{equation*}
G= \frac{x\, (1+x^2)^2}{(1+x+x^2)^3}\ .
\end{equation*}
Note the that the condition $x\neq 1$ is satisfied for $0\leq G< 4/27$. All the expressions above are invariant under $x\to 1/x$, 
so we may always choose $x$ such that $0\leq x<1$.

As a final remark, we note that $t_\infty\equiv \lim_{k\to \infty} t_k=x/(1+x^2)=C$ (the quantity $t_\infty$ enumerates simple slices
with arbitrary left-boundary length $\ell\geq 2$). Recall that the parameter $C$ is nothing but the particular value $t(G)$ of $t$ used in Section~\ref{sec:kernel}
to determine $\tilde{h}_4$ as a function of $G$. The line used  in Section~\ref{sec:kernel} is thus in fact the line $t=t_\infty(G)$. The value $\tilde{\Phi}(t(G))=t(G)/(1+t(G))^2$ that we found may then be undestood as a direct consequence of our recursion relation.
Indeed, letting $k\to \infty$ in our recursion, we may write $t_\infty(G)=(t_\infty(G)+1)\tilde{\Phi}(t_\infty(G))/(1-(t_\infty(G)+1)\tilde{\Phi}(t_\infty(G)))$, which, by inversion, 
reproduces the above value of  $\tilde{\Phi}(t(G))$ when $t(G)=t_\infty(G)$.

\subsection{Solution for general quadrangulations}
\label{sec:general}
Recall the correspondence of eqs.~\eqref{eq:rRtT} and \eqref{eq:gG} between simple and general quadrangulations:
\begin{equation*}
R_k(g)=R_1\, r_k(G)\ , \qquad G=g\, R_1^2
\end{equation*}
and the expression that we just found for $r_k$:
\begin{equation*}
r_k=r_\infty \frac{(1-x^{k})(1-x^{k+3})}{(1-x^{k+1})(1-x^{k+2})}\ , \qquad r_\infty\equiv\frac{1+x+x^2}{1+x^2}\ .
\end{equation*}
with $0\leq x<1$. 
This immediately leads to the expression 
\begin{equation}
R_k=R_\infty \frac{(1-x^{k})(1-x^{k+3})}{(1-x^{k+1})(1-x^{k+2})}\ , \qquad R_\infty\equiv R_1\, r_\infty\  .
\label{eq:finalRk}
\end{equation}
In particular, the prefactor $R_\infty$ being such that $R_\infty=\lim_{k\to \infty} R_k$, it may be interpreted as the generating function of slices 
with arbitrary left-boundary length $\ell\geq 1$. The above definition of $R_\infty$ (via $R_\infty\equiv R_1\, r_\infty$) 
therefore matches precisely that of Section~\ref{sec:GkRk}, hence our notation.
Note also the relation $g\, R_\infty^2= G\, r_\infty^2$. 

As explained in Section~\ref{sec:GkRk}, the value of $R_\infty$ may be obtained directly as the solution of 
the quadratic equation \eqref{eq:Rinf} (satisfying $R_\infty=1+O(g)$), namely
\begin{equation}
R_\infty=\frac{1-\sqrt{1-12g}}{6g}\ .
\label{eq:Rg}
\end{equation}
To fully express $R_k$ in terms of $g$, it remains to connect the parameter $x$ in \eqref{eq:finalRk}
to $g$. Recall the formulas
\begin{equation*}
G=\frac{C}{(1+C)^3}\, \qquad C=\frac{x}{1+x^2}\, \qquad  r_\infty=\frac{1+x+x^2}{1+x^2}=1+C
\end{equation*}
and therefore
\begin{equation*}
g\, R_\infty^2 =G\, r_\infty^2=\frac{C}{1+C}=\frac{x}{1+x+x^2}
\end{equation*}
or equivalently
\begin{equation}
x+\frac{1}{x}+1=\frac{1}{g\, R_\infty^2}
\label{eq:xg}
\end{equation}
with $R_\infty$ as in \eqref{eq:Rg}.
Putting \eqref{eq:Rg} and \eqref{eq:xg} together, we arrive at the following parametrization
\begin{equation}
R_\infty=\frac{1+4x+x^2}{1+x+x^2}\, \qquad g=\frac{x(1+x+x^2)}{(1+4x+x^2)^2}\ .
\label{eq:paramRg}
\end{equation}
The expressions \eqref{eq:finalRk} and \eqref{eq:paramRg} precisely reproduce the result of \cite{GEOD} for $R_k$.
Again demanding $0\leq x<1$ amounts to demanding $0\leq g<1/12$, in agreement with the fact that the
number of quadrangulations with $F$ faces growths exponentially like $12^F$.

From \eqref{eq:GkRk}, we obtain our final formula for the distance-dependent two-point function
\begin{equation*}
\begin{split}
&G_k=\frac{(1-x)^3(1+x)^2 (1+4x+x^2)\, x^{k-1}(1-x^{2k+3})}{(1+x+x^2)(1-x^k)(1-x^{k+1})(1-x^{k+2})(1-x^{k+3})}-\delta_{k,1} \\
&{\rm with}\ \  g=\frac{x(1+x+x^2)}{(1+4x+x^2)^2}\ . \\
\end{split}
\end{equation*}

\section{Conclusion}
\label{sec:conclusion}
To conclude, we notice, as in \cite{G15}, that the decomposition of a slice enumerated by $T_k$ may itself be repeated
recursively inside the sub-slices enumerated by $T_{k-1}$ and so on. This produces (by concatenation of the dividing lines
of the same ``level") a number of nested lines joining the two boundaries of the slice, each line visiting a succession of vertices alternatively at distance
$\ell'-2$ and $\ell'-1$ from the apex for some $\ell'$ ranging from $3$ to the left-boundary length $\ell$ (supposedly being at least $3$)
of the slice at hand. These ``concentric" lines may be viewed as boundaries of the successive balls centered around the apex and 
with radius $\ell'-2$ between $1$ and $\ell-2$ (for some appropriate definition of the balls). More precisely, each ball has also 
in general several \emph{closed} boundaries within the slice encircling connected domains whose vertices are at distance larger than the radius
of the ball. Each of the concentric lines corresponds therefore to a particular boundary, that which separates the apex from the root-vertex of the slice. 
If we complete the ball of radius $\ell'-2$ by the interiors of its closed boundaries, we obtain what can be called the \emph{hull} of radius $\ell'-2$
of the slice. The concentric lines are thus hull boundaries and the statistics of their lengths may in principle be studied
by our formalism\footnote{Note that by closing slices as in figure~\eqref{fig:twopoint}, balls and hulls for slices also correspond to balls and hulls for planar quadrangulations.}.  

Finally, since our approach by recursion was successful in the case of both triangulations and quadrangulations, we may hope that 
a similar scheme could be applied to more general families of maps, for instance maps with prescribed face degrees.



\bibliographystyle{plain}
\bibliography{quadrang2p}

\end{document}